\documentclass[12pt]{article}

\usepackage{amsmath,amssymb,amsfonts,theorem,makeidx,latexsym,epsfig,subfigure,graphics,caption}

\newtheorem{defn}{Definition}[section]

\newtheorem{lemma}[defn]{Lemma}

{\theorembodyfont{\rmfamily}

\newtheorem{ex}[defn]{Example}}

\newtheorem{thm}[defn]{Theorem}

\newtheorem{prop}[defn]{Proposition}

\newtheorem{cor}[defn]{Corollary}

\numberwithin{equation}{section}

\newcommand{\ltr}{ L^2(\mathbb R) }

\newcommand{\mn}{\mathbb N}

\newcommand{\mr}{\mathbb R}

\newcommand{\mz}{\mathbb Z}

\newcommand{\mc}{\mathbb C}

\newcommand{\mt}{ E_{mb}T_{na}g}

\def\bp{{\noindent\bf Proof. \ }}

\def\ep{\hfill$\square$\par\bigskip}

\def\bqs{\begin{equation}}

\def\eqs{\tag*{$\square$}\end{equation}\par\bigskip}

\def\hg{\hat{g}}

\def\la{\langle}

\def\ra{\rangle}

\def\ga{\gamma}

\def\supp{\text{supp}}

\def\bop{\begin{op}\rm}

\def\eop{\end{op}}

\def\bee{\begin{eqnarray}}

\def\ene{\end{eqnarray}}

\def\bes{\begin{eqnarray*}}

\def\ens{\end{eqnarray*}}

\def\bei{\begin{itemize}}

\def\eni{\end{itemize}}

\def\bt{\begin{thm}}

\def\et{\end{thm}}

\def\bc{\begin{cor}}

\def\ec{\end{cor}}

\def\bpr{\begin{prop}}

\def\epr{\end{prop}}

\def\bl{\begin{lemma}}

\def\el{\end{lemma}}

\def\bd{\begin{defn}}

\def\ed{\end{defn}}

\def\bex{\begin{ex}}

\def\enx{\end{ex}}

\def\bfi{\begin{fig}}

\def\efi{\end{fig}}

\def\mt{{\mathbb T}}

\def\hg{\widehat{G}}
\def\ltg{L^2(G)}
\def\lthg{L^2(\hg)}
\def\inhg{\int_{\hg}}
\def\ing{\int_G}

\def\lh{\Lambda^{\bot}}
\def\sulh{\displaystyle{\sum_{\omega \in \lh}}}

\def\lhk{\Lambda_k^\bot}

\def\sulhk{\displaystyle{\sum_{\omega \in \lhk}}}

\def\sulk{\displaystyle{\sum_{\lambda\in \Lambda_k}}}

\def\muhg{\mu_{\hg}}
\def\ml{{\cal M}_\lambda}
\newcommand{\mrs}{{\mathbb R}^s}
\newcommand{\mzs}{{\mathbb Z}^s}
\def\ltrs{L^2(\mrs)}

\def\f{{\cal F}}

\def\hpo{\hat{\psi_0}}

\def\ok{\Omega_k}
\def\oko{\Omega_{k+1}}
\def\vk{V_k}
\def\vko{V_{k+1}}

\def\hko{H_{k+1}}

\def\lak{\Lambda_k}
\def\lako{\Lambda_{k+1}}

\def\mko{M_{k+1}}
\def\nk{N_k}

\title{The unitary extension principle on locally compact abelian groups}

\date{}

\author{Ole Christensen and Say Song Goh}

\begin{document}

\maketitle

\begin{abstract} The unitary extension principle (UEP) by Ron and Shen yields conditions for the construction of a multi-generated tight wavelet frame for $L^2(\mr^s)$ based on a given refinable function.
In this paper we show that the UEP can be generalized to locally compact abelian groups. In the general setting, the resulting frames are generated by modulates of a collection of functions;
via the Fourier transform this corresponds to a generalized shift-invariant system.
Both the stationary and the nonstationary case are covered. We provide general constructions, based on B-splines
on the group itself as well as on characteristic functions on the dual group. Finally, we consider a number of concrete groups and derive explicit constructions of the resulting frames.

\end{abstract}

\section{Introduction} The unitary extension principle (UEP) by Ron and Shen
\cite{RoSh2} and its many variants (\cite{CHS,DRoSh3}, to name a few) are key results in wavelet analysis. They allow construction
of tight wavelet frames with compact support, desired smoothness, and good approximation theoretic properties. In this paper we show how
the theory can be generalized to the setting of locally compact abelian (LCA) groups. The advantage of this approach is threefold. First, we are able to cover several variants of the
unitary extension principle, e.g., the standard case on $\mrs$ and the periodic case corresponding to the group $\mt^s,$ all at once.
Secondly, we can now apply the UEP to a large set of other groups; among these, the group of integers $\mz$ is particularly interesting. Thirdly, the general approach throws new light on the classical UEP:
it reveals the structure that is necessary if we want to consider the UEP from a more general point of view -- a structure that in the case of the group $\mr^s$
turns out to coincide with the classical wavelet structure.

The general approach presented here is able to handle both the stationary and the nonstationary case. In the full generality
of LCA groups we will derive explicit conditions for the UEP construction of tight frames, based on either
B-splines on the group or characteristic functions on the dual group.  Finally, a number of explicit
constructions on groups of particular interest are provided.

Denoting the underlying LCA group by $G,$ the construction yields
frames for $\lthg,$ where $\hg$ denotes the dual group. The frames are obtained by letting a class of modulation operators act on a family of functions in $\lthg.$
Via the Fourier transform, this immediately yields a generalized shift-invariant system
that forms a frame for $L^2(G);$ see the more technical description right after the definition of the modulation operator in \eqref{11h}.

The discussion of the general results is complemented by  explicit constructions
on the LCA groups $\mrs, \mt, \mz,$ and $\mz_N$ (the integers modulo $N$). Some of them are based on a generalization of B-splines to LCA groups
that was discovered already in 1994 (independently by Dahlke \cite{Dah} and Tikhomirov \cite{T}); other constructions  are generated by characteristic functions for certain sets in $\hg.$

Section \ref{17a} will provide us with the necessary background on LCA groups; the general version of the UEP and its proof is stated in Section \ref{2101a}.
The formulation of the results based on B-splines  are presented in Section \ref{1901d}, while the case of characteristic functions on $\hg$ are
in Section \ref{471225a}; in both cases  applications to a number of specific LCA groups
are given.

Let us end this introduction with a few comments about technicalities. The main
difficulty in the extension of the UEP to LCA groups is that there is no
scaling operator on LCA groups. The way to overcome this issue turns out to be to
consider a collection of modulation operators acting on a family of frame generators,
rather than the usual collection of scaled and translated versions of a fixed function. This leads to a different form of the scaling equation. Also, the
traditional assumption of the wavelet subspaces being nested has to be replaced by
a condition on a nested sequence of lattices in the group. After getting familiar
with these new aspects, the reader will observe that several of the technical results
follow the same pattern as in the classical proofs of the UEP.

We also note that the approach in the paper heavily uses the LCA group structure.
We would like to point the attention of the reader to a different generalization
of the UEP, taking place on smooth and compact Riemannian manifolds; see \cite{WZ}.
We  also mention that there is a growing literature on wavelet analysis on the
$p$-adic numbers; see \cite{Skopina}
and the references therein.
The $p$-adic numbers do not have nontrivial
lattices \cite{BowRos2014}, and are not covered by our methods.  A more general wavelet theory on
LCA groups with a compact open subgroup is developed in \cite{benedetto}.

\section{Preliminaries on LCA groups} \label{17a}
In this section we will give a short introduction to the necessary background on LCA groups; for more information we refer to the books \cite{HeRo,ReSt,Ru3}.

Let $G$ be an LCA group, with the group
composition denoted by the symbol $``+"$ and the neutral element 0. We will assume that $G$ is equipped with a Hausdorff topology, and that $G$
is a countable union of compact sets and metrizable. A {\it character} on $G$ is a function
$\ga: G \to \mt:=\{z\in \mc \ \big| \ |z|=1\}$ that satisfies the condition
$\ga(x+y)=\ga(x) \ga(y),\ x,y\in G.$ We denote the set of continuous characters by $\hg,$ which also forms an LCA group, the {\it dual group} of $G,$  when equipped with the composition
$(\ga + \ga^\prime)(x):= \ga(x) \ga^\prime(x),$ $\ga, \ga^\prime\in \hg,$ $x\in G,$ and an appropriate  topology.
It is a classical result that the double-dual group $\widehat{\hg}$ is topologically isomorphic to the group $G;$ usually we can identify the two groups and we will simply write
$\widehat{\hg}=G.$ Thus  $\ga(x)$ can be interpreted as either the action of $\ga\in \hg$ on $x\in G,$ or the action of
$x\in \widehat{\hg}=G$ on $\ga \in \hg;$ for this reason we will from now on use the notation
\bes (x,\ga):= \ga(x), \ x\in G, \ \ga\in \hg.\ens

The LCA group $G$ can be equipped with a Radon measure $\mu_G$ that is {\it translation invariant}, which means that for all
continuous functions $f$ on $G$ with compact support,
\bee \label{11e} \ing f(x+y)\, d\mu_G(x)=\ing f(x)\, d\mu_G(x), \ \forall \, y\in G.\ene This measure is unique up to scalar multiplication, and is called the
{\it Haar measure.}
We will consider the Haar measure $\mu_G$ as fixed throughout the paper.
Based on the Haar measure, we define the spaces
$L^1(G), L^2(G)$ and $L^\infty(G)$ in the usual way.
The space $\ltg$ is a Hilbert space, and our assumption of $G$ being a countable union of compact sets and metrizable implies (and is, in fact, equivalent to) $\ltg$ being separable.

The {\it Fourier transform} is defined by
\bee
\label{2.1aaa}
 {\cal F}: L^1(G) \to C_0(\hg), \quad {\cal F}f(\ga):= \hat{f}(\ga) := \ing f(x)(-x, \ga)\,d\mu_G(x).
\ene
The {\it inversion theorem} states that with appropriate normalization of the Haar measure $\muhg$ on $\hg,$ for $f\in L^1(G)$ such that $\hat{f}\in L^1(\hg),$ it holds that
\bee \label{11f} f(x)= \inhg \hat{f}(\ga) \, (x,\ga) d\muhg(\ga), \ x\in G.\ene
We will always normalize the measure on $\hg$ in such a way that the inversion formula holds. With this choice,
the Fourier transform  can be extended to a surjective isometry ${\cal F}: L^2(G) \to L^2(\hg),$  exactly as in the classical case of $G=\mr.$
To simplify notations, from now onwards, in all integrals when the context is clear (e.g., (\ref{11e})--(\ref{11f})), we simply write $d\mu_G(x) = dx$ and $d\mu_{\widehat{G}}(\ga) = d\ga$.

Among the examples of LCA groups, we find $\mr,$ $\mt,$ $\mz,$ $\mz_N$,
as well as their higher-dimensional variants and direct products hereof; following \cite{ReSt} and \cite{FeiKo} we will call such groups {\it elementary LCA groups.} As discussed in \cite{FeiKo},
various typical problems in signal processing can be modeled using elementary LCA groups.

A {\it lattice} (sometimes called a uniform lattice) in an LCA group $G$ is a discrete subgroup $\Lambda$
for which $G/\Lambda$ is compact. Lattices are known explicitly in all the elementary LCA groups  and in many other  LCA groups; however,  there also exist  LCA groups without lattices, see, e.g., \cite{KK,BowRos2014,LeSi1}.
The {\it annihilator} $\Lambda^\bot$ of a lattice $\Lambda$  is defined by
\bes \Lambda^\bot:= \{ \ga \in \hg \ \big| \ (x, \ga)=1, \ \forall x\in \Lambda\}.\ens
It follows from the definition of the topology on $\hg$ that the annihilator $\Lambda^\bot$ is a closed subgroup of $\hg.$
A lattice in $G$ leads to a splitting of the groups $G$
and $\hg$ into disjoint cosets,  see, e.g., \cite{KK}:

\bl \label{11c} Let $G$ be an LCA group and $\Lambda$ a lattice in $G.$ Then the following hold:

\bei \item[{\rm (i)}] There exists a Borel measurable relatively compact set $Q\subseteq G$ such that
\bee \label{11a} G= \bigcup_{\lambda\in \Lambda} (\lambda+Q), \ \ (\lambda+Q) \cap (\lambda^\prime +Q) = \emptyset \ \mbox{for} \ \lambda\neq \lambda^\prime, \ \lambda, \lambda^\prime \in \Lambda.\ene
\item[{\rm (ii)}] The set $\lh$ is a lattice in $\hg,$ and there exists a Borel measurable relatively compact set $V\subseteq \hg$ such that
\bes
\hg= \bigcup_{\omega\in \lh} (\omega+V), \ \ (\omega+V)  \cap (\omega^\prime +V) = \emptyset \ \mbox{for} \ \omega\neq \omega^\prime,
\ \omega, \omega^\prime \in \lh.\ens \eni\el

A set $Q$ as in \eqref{11a} that has the properties in Lemma \ref{11c}(i) is called a {\it fundamental domain} associated to the lattice $\Lambda.$ For convenience
we will allow sets $Q$ for which the two conditions in \eqref{11a} hold up to a set of measure zero.
When we speak about a {\it periodic function} $f\in L^\infty(Q),$ it is understood that we extend $f$ to a function on $G$ by
\bes
f(\lambda + x):= f(x), \, \lambda \in \Lambda, \, x\in Q.\ens

Given a lattice $\Lambda$ in $G,$ choose a fundamental domain
$Q.$ The {\it density} of $\Lambda$ is defined by
$s(\Lambda):= \mu_G(Q).$ It is well known that this is independent of the chosen fundamental domain, and that (see \cite{G6}) $s(\Lambda) \, s(\Lambda^\bot)=1.$

Given any $\lambda\in G,$ consider the
{\it generalized modulation operator}
\bee \label{11h} \ml: \lthg \to \lthg, \ (\ml f)(\ga):= (\lambda, \ga) \, f(\ga).\ene As for the
modulation operator on $\mr,$ it is easy to see that
$\ml$ is a unitary operator.

The main outcome of the current paper is a method for constructing frames
for $\lthg$ of the form $\{\ml \Gamma_k\}_{k\in J, \lambda\in \Lambda_k},$ where
$\{\Gamma_k\}_{k\in J}$ is a countable collection of functions in $\lthg$ and
$\{\Lambda_k\}_{k\in J}$ a collection of lattices in $G.$ Considering the {\it translation operator} $T_y,$ $y \in G,$ on $L^2(G)$ defined by
$ T_y: L^2(G) \to L^2(G), \, T_yf(x):=f(x-y), \, x\in G,$
it is well known that
$ \label{99161a}
\f T_y= {\cal M}_{-y} \f;$
thus, if $\{\ml \Gamma_k\}_{k\in J, \lambda\in \Lambda_k}$ is a frame for $L^2(\widehat{G})$
it immediately follows that the system
$ \{{\cal F}^{-1}\ml \Gamma_k\}_{k\in J, \lambda\in \Lambda_k}=
\{T_{\lambda} {\cal F}^{-1}\Gamma_k\}_{k\in J, \lambda\in \Lambda_k}$ is a frame for
$L^2(G).$ This system is a so-called {\it generalized shift-invariant system}
in the terminology of Ron and Shen \cite{RoSh7}. The case where the lattices
$\Lambda_k$ are independent of $k$ correspond to the classical shift-invariant
systems; for a detailed analysis of such systems we refer to \cite{CP}.

We will need the following result, which is a variant of Lemma 3.3 in \cite{KL}; the
version stated here is proved in \cite{CG} and repeated in \cite{CBN}. For technical reasons in the lemma and subsequently in Section \ref{2101a},
we consider the dense subspace
$C_c(\hg)$ of $\lthg$ defined by
\bes
C_c(\hg):= \{f\in L^2(\widehat{G}) \, \big| \, f \, \mbox{is continuous and compactly supported}\}.\ens

\bl \label{11m} Let $\Lambda$ be a lattice in $G,$ and let $V \subseteq \hg$ denote a fundamental domain associated with the lattice $\lh.$ Let $F, \Phi \in \lthg.$ Then the function
\bes
\alpha : \hg \to \mc, \ \alpha(\ga):= \sulh F(\omega + \ga) \overline{\Phi(\omega + \ga)}\ens is well defined,  belongs  to
$L^1(V),$ and satisfies that
\bes
\alpha(\ga+ \omega^\prime)= \alpha(\ga), \ \forall \, \ga \in \hg,  \omega^\prime \in \lh.\ens In addition, if $F \in C_c(\hg),$ then
\bes \sum_{\lambda \in \Lambda} |\la F, \ml \Phi \ra|^2 = \mu_{\hg}(V)\,  \int_V \bigg|
\sulh F(\omega + \ga) \overline{\Phi(\omega + \ga)}\bigg|^2 \, d\ga. \ens
 \el

One of the central ingredient in our generalization of the UEP is to consider
a family of nested lattices $\{\Lambda_k\}_{k\in I}$ in the group $G.$ We will need
to consider relations between the lattices and their corresponding fundamental domains;
for easy reference we will formulate here the relevant connections between just two lattices.

\bl \label{801b} Consider two lattices $\Lambda_0 \subset \Lambda_1$ in an LCA group $G.$ Then the following hold:

\bei
\item[{\rm (i)}] The quotient group
$\Lambda_1/\Lambda_0$ is finite, with cardinality
$ \label{1501c} \left| \Lambda_1/\Lambda_0\right| = \frac{s(\Lambda_0)}{s(\Lambda_1)}.$
\item[{\rm (ii)}] Let $\{\eta_\ell\}_{\ell=1}^d\subset \Lambda_1$ be a fundamental domain for
$\Lambda_0,$ considered as a subgroup of $\Lambda_1,$ chosen such that
    $\eta_1=0.$ Then
     \bes
     \Lambda_1= \bigcup_{\ell=1}^d (\eta_\ell + \Lambda_0),
    \ \mbox{with} \ (\eta_\ell + \Lambda_0) \cap (\eta_{\ell^\prime} + \Lambda_0)=
    \emptyset \ \mbox{if} \ \ell\neq \ell^\prime.\ens
\item[{\rm (iii)}]Choose $Q_1$ as a fundamental domain associated with the lattice
$\Lambda_1$ in $G.$ Then \bee \label{801da} (\eta_\ell + Q_1) \cap (\eta_{\ell^\prime} + Q_1)= \emptyset \ \mbox{if} \ \ell\neq \ell^\prime.\ene Furthermore, the set
$Q_0:= \bigcup_{\ell=1}^d (\eta_\ell + Q_1)$ is a fundamental domain associated with the lattice $\Lambda_0$ in $G.$
\item[{\rm (iv)}]
$\Lambda_1^\bot$ is a subgroup of $\Lambda_0^\bot;$ the quotient group
$\Lambda_0^\bot/\Lambda_1^\bot$ is finite, and
$ \label{1501b} \left| \Lambda_0^\bot/\Lambda_1^\bot\right|=\left| \Lambda_1/\Lambda_0\right|.$
\item[{\rm (v)}]  Let $\{\nu_\ell\}_{\ell=1}^d\subset \Lambda_0^\bot$ be a fundamental domain for $\Lambda_1^\bot,$ considered as a subgroup of $\Lambda_0^\bot,$ chosen such that
    $\nu_1=0.$ Then
    \bes
    \Lambda_0^\bot= \bigcup_{\ell=1}^d (\nu_\ell + \Lambda_1^\bot),
    \ \mbox{with} \ (\nu_\ell + \Lambda_1^\bot) \cap (\nu_{\ell^\prime} + \Lambda_1^\bot)=
    \emptyset \ \mbox{if} \ \ell\neq \ell^\prime.\ens
\item[{\rm (vi)}] Let $V_0\subset \hg$ be a fundamental domain associated with the lattice $\Lambda_0^\bot$ in $\hg.$ Then
\bes
(\nu_\ell + V_0) \cap (\nu_{\ell^\prime} + V_0)= \emptyset \ \mbox{if} \ \ell\neq \ell^\prime.\ens
Furthermore, the set $V_1:= \bigcup_{\ell=1}^d (\nu_\ell + V_0)$ is a fundamental domain associated with the lattice $\Lambda_1^\bot$ in $\hg.$
\eni \el

\bp (i) By definition of a lattice we know that $G/\Lambda_0$ is relatively compact; since
$\Lambda_1\subseteq G$ this implies that $\Lambda_1/\Lambda_0$ is also relatively compact. But since $\Lambda_1$ is discrete, $\Lambda_1/\Lambda_0$ is also discrete;
hence the set must be finite. We postpone the proof of the cardinality by a few lines.

\vspace{.1in}\noindent(ii) By definition of the fundamental domain,
$ \Lambda_1= \bigcup_{\lambda\in \Lambda_0}  (\lambda + \{\eta_\ell\}_{\ell=1}^d) =
\bigcup_{\ell=1}^d (\eta_\ell + \Lambda_0).$ By construction the union is disjoint, so (ii) holds.

\vspace{.1in}\noindent(iii) As $Q_1$ is a fundamental domain associated with
$\Lambda_1,$ we know that $(\lambda + Q_1) \cap (\lambda^\prime + Q_1)= \emptyset$ if
$\lambda, \lambda^\prime \in \Lambda_1, \, \lambda \neq \lambda^\prime.$
Since $\eta_\ell \in \Lambda_1,$ \eqref{801da} follows.
In order to prove the rest of (iii), using (ii) we see that
\bes G = \bigcup_{\lambda \in \Lambda_1} (\lambda + Q_1)=
\bigcup_{\lambda \in \Lambda_0} \bigcup_{\ell=1}^d (\eta_\ell+\lambda + Q_1)
 & = & \bigcup_{\lambda \in \Lambda_0} \big(\lambda + \bigcup_{\ell=1}^d (\eta_\ell + Q_1)\big)= \bigcup_{\lambda \in \Lambda_0} (\lambda + Q_0).\ens
 Now assume that
 for some $\lambda, \lambda^\prime \in \Lambda_0$ we have
 $\mu_G\big( (\lambda + Q_0) \cap (\lambda^\prime + Q_0) \big) > 0.$ Then for some
 $\eta_\ell, \eta_{\ell^\prime}\in \Lambda_1$ we have
 $\mu_G \big( (\lambda + \eta_\ell + Q_1) \cap (\lambda^\prime +\eta_{\ell^\prime}+ Q_1) \big) > 0.$ Since $\lambda + \eta_\ell, \lambda^\prime +\eta_{\ell^\prime}\in \Lambda_1,$
 we conclude that $\lambda + \eta_\ell= \lambda^\prime +\eta_{\ell^\prime};$
 by (ii), this implies  that $\eta_\ell= \eta_{\ell^\prime}$  and therefore  $\lambda= \lambda^\prime,$ i.e., $Q_0$ is indeed a
 fundamental domain associated with $\Lambda_0.$

Let us now give the proof of the cardinality in (i). Note that the cardinality
of the set $\Lambda_1/\Lambda_0$ equals the number $d$ introduced in (ii). Since $Q_1$ is a fundamental domain associated with $\Lambda_1,$ by definition we have
$s(\Lambda_1)= \mu_G(Q_1).$ It now follows from (iii), the disjointness (up to a set of measure zero)
of the sets $\eta_\ell+Q_1,$ $\ell = 1, \ldots, d,$ and the translation invariance of the measure  that
$$ s(\Lambda_0)= \mu_G (Q_0)= \mu_G\big(\bigcup_{\ell=1}^d (\eta_\ell + Q_1)\big) = d\, \mu_G(Q_1)
= d\, s(\Lambda_1),$$ as claimed.

\vspace{.1in}\noindent(iv) This follows from Proposition 4.2.24 in \cite{ReSt}, but let us give a direct proof.
By the definition of the annihilator, the assumption
$\Lambda_0 \subset \Lambda_1$ implies that $\Lambda_1^\bot \subset \Lambda_0^\bot;$
since both are groups it is clear that $\Lambda_1^\bot$ is a subgroup of $\Lambda_0^\bot.$ That $\Lambda_0^\bot/\Lambda_1^\bot$ is finite now follows from (i).
Using (i), we arrive at
$\left| \Lambda_0^\bot/\Lambda_1^\bot\right| = \frac{s(\Lambda_1^\bot)}{s(\Lambda_0^\bot)}=\frac{s(\Lambda_0)}{s(\Lambda_1)}
=\left| \Lambda_1/\Lambda_0\right|,$ which proves (iv).

\vspace{.1in}\noindent(v), (vi)
Finally, the results in (v) and (vi) follow immediately from (ii) and (iii) applied to the inclusion $\Lambda_1^\bot \subset \Lambda_0^\bot.$ \ep

\section{The unitary extension principle} \label{2101a}

In order to avoid a long list of assumptions in the formulation of the UEP, we will state  the standing assumptions for this section in  a ``General setup".
Before we do this, let us mention a few conventions that will help us avoid cumbersome notations.

First, the UEP will be based on a sequence of functions
$\{ \Phi_k\}_{k\in I},$ indexed by
a countable sequence of consecutive integers in $\mz,$ i.e., either
$I= \{k\}_{k= k_0}^\infty$ or  $I= \{k\}_{k= k_0}^{k_1}$ for some $k_0, k_1\in \mz.$
In what follows we will tacitly assume that  $I=\{k\}_{k= k_0}^\infty$ and leave the minor modifications in the case $I= \{k\}_{k= k_0}^{k_1}$ to the reader. A typical example: a
condition involving $\Phi_k$ and $\Phi_{k+1}$ makes perfect sense for
$k\in I$ if we assume that $I= \{k\}_{k= k_0}^\infty.$ On the other hand, for the case
$I= \{k\}_{k= k_0}^{k_1}$ one would have to assume that $k\in \{k_0, \dots, k_1-1\}.$
For the rest of the paper, we will let $k_0$ denote the starting index of the set $I.$

The UEP on $\mrs$ by Ron and Shen is formulated in terms of conditions on some filters, which are periodic functions. We will need an analog concept in our setting, and we will use the convention stated right after Lemma \ref{11c}. As we will see, the relevant periodic functions are actually defined on the dual group $\hg.$

Let us now state the standing assumptions for this section:

\vspace{.1in}\noindent{\bf General setup:}  Let
$I$ be a sequence of consecutive numbers in $\mz.$
Let
$\{ \Lambda_k\}_{k\in I}$ be a nested sequence of lattices in $G,$ i.e.,
\bee \label{901k} \Lambda_{k_0} \subset \Lambda_{k_0+1} \subset \Lambda_{k_0+2} \subset \cdots.\ene
Let $\{ \Phi_k\}_{k\in I}$ be a sequence of functions in $\lthg.$ Furthermore, for each $k\in I,$
let $V_k$ denote a fundamental domain associated with the lattice $\lhk;$ then, in particular, for each $k\in I,$
\bee \label{901f} \hg = \bigcup_{\omega\in \lhk} (\omega + V_k),
\ (\omega+V_k)  \cap (\omega^\prime +V_k) = \emptyset \ \mbox{for} \ \omega\neq \omega^\prime, \ \omega, \omega^\prime \in \lhk.\ene

Assume the following conditions:
\bei
\item[(i)] For every compact set $S$ in $\hg,$ there exists $K_1\in I$ such that
    \bee \label{901g}
    \mu_{\widehat{G}} \big( (\omega + S) \cap (\omega^\prime + S) \big) = 0 \ \mbox{for} \ \omega \neq \omega^\prime, \ \omega, \omega^\prime \in \Lambda_{K_1}^\bot.
    \ene
\item[(ii)]  For every compact set $S$ in $\hg$ and any $\epsilon >0,$ there exists $K_2\in I$ such that for all $k\ge K_2,  \, k\in I,$
    \bee \label{901a} \big| \mu_{\hg}(V_k) \, |\Phi_k(\ga) |^2 -1 \big| \le \epsilon, \ \forall \ga\in S.\ene
\item[(iii)] For all $k \in I$ and some periodic functions $H_{k+1} \in L^\infty(V_{k+1})$ (see the convention in Section \ref{17a}),
    \bee \label{1001a} \Phi_{k}(\ga)= H_{k+1}(\ga) \, \Phi_{k+1}(\ga), \, \mbox{a.e.} \,
        \ga \in \hg.\ene
\eni

For $k\in I,$ given periodic functions $G_{k+1}^{(m)}\in L^\infty(V_{k+1}), \, m=1, \dots, \rho_k,$ define the
functions $\Psi_k^{(m)} \in \lthg, \, m=1, \dots, \rho_k,$ by
\bee \label{1001b} \Psi_k^{(m)}(\ga):= G_{k+1}^{(m)}(\ga)\,  \Phi_{k+1}(\ga), \ \ga \in \hg.\ene
Our goal is to identify conditions on the {\it filters} $H_k$ and $G_k^{(m)}$ such that
the collection of functions
\bes
\{\ml \Phi_{k_0}\}_{\lambda\in \Lambda_{k_0}} \bigcup \,
\{\ml \Psi_k^{(m)}\}_{k\ge k_0, \, \lambda\in \Lambda_{k}, m=1, \dots, \rho_k}\ens forms a tight frame for $\lthg$ with frame bound 1.

\vspace{.3in} We will see later in Example \ref{eg3.7} that for the case $G=\mr,$ the assumptions \eqref{901a} and \eqref{1001a}
correspond directly to the assumptions and setup used in the classical UEP constructions; on the other hand, the condition \eqref{901g} is automatically satisfied in this case,
and does not appear explicitly in the classical UEP.

Starting with a suitable choice of a sequence
$\{\nu_{k_0, \ell}\}_{\ell=1, \dots, d_{k_0}}$ associated with the
``lowest'' level lattice $\Lambda_{k_0}^\bot,$ repeated use of Lemma \ref{801b}(v) shows that for each $k\in I,$ we can choose a sequence $\{\nu_{k, \ell}\}_{\ell=1, \dots, d_k} \subset \hg$ such
that $\nu_{k,1}=0$ and
\bee \label{1001d} \Lambda_k^\bot & = & \bigcup_{\ell=1}^{d_k}( \nu_{k, \ell}+ \Lambda_{k+1}^\bot), \ \  ( \nu_{k, \ell}+ \Lambda_{k+1}^\bot) \cap   ( \nu_{k, \ell^\prime}+ \Lambda_{k+1}^\bot) = \emptyset \ \mbox{for} \
\ell \neq \ell^\prime.
\ene
For $k\in I,$ consider the $(\rho_k+1) \times d_k$ matrix-valued function $P_k$ defined by
\bee \label{1001g} P_k(\ga):= \begin{pmatrix} H_{k+1} (\ga + \nu_{k,1}) & \cdot & \cdot & \cdot & H_{k+1} (\ga + \nu_{k,d_k}) \\
G_{k+1}^{(1)} (\ga + \nu_{k, 1}) & \cdot & \cdot & \cdot & G_{k+1}^{(1)} (\ga + \nu_{k,d_k}) \\ \cdot & \cdot & \cdot & \cdot & \cdot \\
 \cdot & \cdot & \cdot & \cdot & \cdot \\
G_{k+1}^{(\rho_k)} (\ga + \nu_{k, 1}) & \cdot & \cdot & \cdot & G_{k+1}^{(\rho_k)} (\ga + \nu_{k,d_k})
\end{pmatrix}, \ \ga \in V_k.\ene

By Lemma \ref{801b}(vi) applied to the lattices $\Lambda_k^\bot$ and $\Lambda_{k+1}^\bot,$ we know that for any fundamental domain $V_k$ associated with $\Lambda_k^\bot,$
\bee \label{172b} (\nu_{k, \ell} + V_k) \, \cap \, (\nu_{k, \ell^\prime} + V_k)= \emptyset
\ \mbox{for} \ \ell \neq \ell^\prime.\ene
Lemma \ref{801b}(vi) also shows that the set
\begin{equation} \label{172c} V_{k+1}^\prime:=  \bigcup_{\ell=1}^{d_k} ( \nu_{k,\ell}+V_k)
\end{equation} is a fundamental domain associated with the lattice $\Lambda_{k+1}^\bot.$ This observation turns out to be important: in fact, some of the analysis to follow applies whenever
$\{V_k\}_{k\in I}$ is an arbitrary collection of fundamental domains associated with the lattices $\{\lak^\bot\}_{k\in I},$ but some of the results require a relationship between fundamental domains ``on consecutive levels'';
in such cases we will apply \eqref{172c}.

Since we will work with different choices of  fundamental domains $V_k,$ we need the following elementary result:

\bl \label{172a} For any two  fundamental domains $V_k$ and $V_k^\prime$ associated with $\Lambda_k^\bot,$ we have $\mu_{\hg}(V_{k})= \mu_{\hg}(V_k^\prime)$ and
$L^\infty(V_k)= L^\infty(V_k^\prime).$\el

Lemma \ref{172a} will help us to simplify the notation: indeed, whenever we need
to deal with, e.g., the fundamental domain $\vko^\prime$ in \eqref{172c}, we will always
write $\mu_{\hg}(V_{k+1})$ instead of  $\mu_{\hg}(V_{k+1}^\prime)$ and
$L^\infty(V_{k+1})$ instead of $L^\infty(V_{k+1}^\prime).$

Let us comment on the fact that $P_k$ for the moment is only considered on $V_k.$
Looking at the entries of the matrix $P_k,$ we see that by considering $\ga \in V_k,$ we exactly
use the information about the functions $H_{k+1}$ and $G_{k+1}^{(m)}$ on the fundamental
domain $V_{k+1}^\prime$ defined by \eqref{172c}; thus, considering $P_k$ only on $V_k$ is consistent with our desire to consider  $H_{k+1}$ and $G_{k+1}^{(m)}$ as periodic functions in  $L^\infty(V_{k+1}).$

We will see that the essential condition in the UEP is that
\bee \label{172e} P_k(\ga)^* P_k(\ga)= d_k I_{d_k}, \, \mbox{a.e.} \, \ga \in V_k.\ene
Note that $d_k=\frac{\mu_{\hg}(V_{k+1})}{\mu_{\hg}(V_k)}.$ We will now show that if \eqref{172e}
holds on {\it any} fundamental domain $V_k,$ then it automatically holds on $\widehat{G}.$

\bl \label{172d} Let $V_k$ denote any fundamental domain associated with
$\Lambda_k^\bot.$ If \eqref{172e} holds for a.e.\ $\ga \in V_k,$ then \eqref{172e} holds for a.e.\ $\ga \in \widehat{G}.$ \el

\bp  The condition \eqref{172e} means that for a.e.\ $\ga \in V_k$ and $\ell, \ell^\prime=1, \dots, d_k,$
\bee \label{182a} H_{k+1}(\ga + \nu_{k,\ell}) \overline{H_{k+1}(\ga + \nu_{k,\ell^\prime})}+
\sum_{m=1}^{\rho_k} G_{k+1}^{(m)}(\ga + \nu_{k,\ell}) \overline{G_{k+1}^{(m)}(\ga + \nu_{k,\ell^\prime})}= d_k \delta_{\ell, \ell^\prime}.\ene
Consider now any $\ga^\prime$ belonging to the fundamental domain $V_{k+1}^\prime$ chosen in
\eqref{172c}. Such  $\ga^\prime$ can be written in the form
$\ga^\prime = \ga +\nu_{k, \tilde{\ell}}$ for some
$\ga \in V_k $ and some $\tilde{\ell} \in \{1, \dots, d_k\}.$ Using \eqref{182a},
that $\nu_{k, \tilde{\ell}}+ \nu_{k, \ell} \in \Lambda_k^\bot= \cup_{q=1}^{d_k}
(\nu_{k, q}+ \Lambda_{k+1}^\bot),$ and the periodicity of
the functions $H_{k+1}$ and $G_{k+1}^{(m)}$ in
$L^\infty(V_{k+1}),$ it follows that $\ga^\prime$ also satisfies \eqref{182a} (we leave the details to the reader).
Using now the periodicity of the entries of the matrix $P_k,$ it follows that \eqref{182a} holds  for a.e.\ $\ga \in \widehat{G}.$
\ep

We now state a lemma which allows us to move around between different levels of
$\Phi_k.$

\bl \label{901bf} In addition to the general setup, assume that for some $k\in I,$ the
matrix-valued function $P_k$ satisfies \eqref{172e}.  Then for all
$F\in C_c(\hg),$
\bee \label{1001m} \sum_{\lambda\in \Lambda_{k+1}} | \la F, \ml \Phi_{k+1}\ra|^2= \sum_{\lambda\in \Lambda_{k}} | \la F, \ml \Phi_k\ra|^2+
\sum_{m=1}^{\rho_k}\sum_{\lambda\in \Lambda_{k}} | \la F, \ml \Psi_k^{(m)}\ra|^2.\ene
\el

\bp Fix $k\in I.$ We first note that by Lemma \ref{172d} the condition \eqref{172e}
holds for a.e.\ $\ga \in \widehat{G}.$ We will now use Lemma \ref{11m} to rewrite the three expressions appearing in \eqref{1001m}. For technical reasons we use the
given fundamental domain $V_k$ for the terms on the right-hand side, and the domain
$V_{k+1}^\prime$ in \eqref{172c} for the left-hand side:
\bee \label{1001ma} \sum_{\lambda\in \Lambda_{k+1}} | \la F, \ml \Phi_{k+1}\ra|^2 & = &
\mu_{\hg}(V_{k+1}) \, \int_{V_{k+1}^\prime} \bigg| \sum_{\omega\in \Lambda_{k+1}^\bot} F(\omega + \ga) \overline{\Phi_{k+1}(\omega + \ga)} \bigg|^2 \, d\ga; \\ \label{1001n}
\sum_{\lambda\in \Lambda_{k}} | \la F, \ml \Phi_{k}\ra|^2 & = &
\mu_{\hg}(V_k) \, \int_{V_k} \bigg| \sum_{\omega\in \Lambda_k^\bot} F(\omega + \ga) \overline{\Phi_{k}(\omega + \ga)} \bigg|^2 \, d\ga; \\ \label{1001p}
\sum_{m=1}^{\rho_k}\sum_{\lambda\in \Lambda_{k}} | \la F, \ml \Psi_{k}^{(m)}\ra|^2 & = &
\sum_{m=1}^{\rho_k}\mu_{\hg}(V_k) \, \int_{V_k} \bigg| \sum_{\omega\in \Lambda_k^\bot} F(\omega + \ga) \overline{\Psi_{k}^{(m)}(\omega + \ga)} \bigg|^2 \, d\ga.
\ene
Based on these expressions, a natural approach is to start with \eqref{1001ma} and
apply the disjoint splitting of $V_{k+1}^\prime$ in terms of $V_k,$ see \eqref{172c}; this
yields
\bee \notag & \ & \sum_{\lambda\in \Lambda_{k+1}} | \la F, \ml \Phi_{k+1}\ra|^2  =
\mu_{\hg}(V_{k+1}) \, \int_{V_{k+1}^\prime} \bigg| \sum_{\omega\in \Lambda_{k+1}^\bot} F(\omega + \ga) \overline{\Phi_{k+1}(\omega + \ga)} \bigg|^2 \, d\ga \\ \notag & = &
\mu_{\hg}(V_{k+1}) \, \sum_{\ell=1}^{d_k}\int_{\nu_{k,\ell}+V_{k}} \bigg| \sum_{\omega\in \Lambda_{k+1}^\bot} F(\omega + \ga) \overline{\Phi_{k+1}(\omega + \ga)} \bigg|^2 \, d\ga \\ \label{1301b} & = &
\mu_{\hg}(V_{k+1}) \, \int_{V_{k}} \sum_{\ell=1}^{d_k}\bigg| \sum_{\omega\in \Lambda_{k+1}^\bot} F(\omega +\nu_{k,\ell} +\ga) \overline{\Phi_{k+1}(\omega + \nu_{k,\ell}+\ga)} \bigg|^2 \, d\ga. \ene

Define the functions $a_k^\ell, \, \ell=1, \dots, d_k,$ and $b_k$ on $V_k$ by
\bes a_k^\ell(\ga) & := & \sum_{\omega\in \Lambda_{k+1}^\bot} F(\omega + \nu_{k, \ell}+ \ga) \overline{\Phi_{k+1}(\omega + \nu_{k, \ell}+ \ga)}, \ \ga \in V_k;
\\ b_k(\ga) & := &  \sum_{\omega\in \Lambda_{k}^\bot} F(\omega + \ga) \overline{\Phi_{k}(\omega + \ga)}, \ \ga \in V_k.\ens
By Lemma \ref{11m}, these functions are well defined.
We will also consider the vector
\bes  c_k(\ga) & := &  \begin{pmatrix} \sum_{\omega\in \Lambda_{k}^\bot} F(\omega + \ga) \overline{\Psi_{k}^{(1)}(\omega + \ga)}  \\
\cdot \\ \cdot \\
\sum_{\omega\in \Lambda_{k}^\bot} F(\omega + \ga) \overline{\Psi_{k}^{(\rho_k)}(\omega + \ga)} \end{pmatrix}, \ \ga \in V_k.\ens
Using \eqref{1001a} and that $H_{k+1}(\ga + \omega) = H_{k+1}(\ga)$ for $\ga \in \hg, \, \omega \in \Lambda_{k+1}^\bot,$
\bes \notag \sum_{\ell=1}^{d_k} a_k^\ell(\ga) \overline{H_{k+1}(\ga + \nu_{k, \ell})}
& = &  \sum_{\ell=1}^{d_k} \sum_{\omega\in \Lambda_{k+1}^\bot} F(\omega + \nu_{k, \ell}+ \ga) \overline{\Phi_{k+1}(\omega + \nu_{k, \ell}+ \ga) H_{k+1}(\omega+\ga + \nu_{k, \ell} )}
\\ 
& = &  \sum_{\ell=1}^{d_k} \sum_{\omega\in \Lambda_{k+1}^\bot} F(\omega + \nu_{k, \ell}+ \ga) \overline{\Phi_{k}(\omega + \nu_{k, \ell}+ \ga)  }.\ens
Via the disjoint splitting of $\Lambda_k^\bot$ in \eqref{1001d}, it follows that
\bes \sum_{\ell=1}^{d_k} a_k^\ell(\ga) \overline{H_{k+1}(\ga + \nu_{k, \ell})}
& = & \sum_{\omega\in \Lambda_{k}^\bot} F(\omega + \ga) \overline{\Phi_{k}(\omega + \ga)  }=b_k(\ga).\ens
In the same way it can be proved from (\ref{1001b}) that
\bes
\sum_{\ell=1}^{d_k} a_k^\ell(\ga)
\begin{pmatrix} \overline{G_{k+1}^{(1)}(\ga + \nu_{k, \ell})} \\ \cdot \\ \cdot \\
\overline{G_{k+1}^{(\rho_k)}(\ga + \nu_{k, \ell})}\end{pmatrix}=c_k(\ga).\ens
Defining the column vectors
\bes \beta_k(\ga):= \begin{pmatrix} b_k(\ga) \\ c_k(\ga) \end{pmatrix}, \ \
\alpha_k(\ga):=  \begin{pmatrix} a_k^1(\ga) \\ \cdot \\ \cdot \\ a_k^{d_k}(\ga) \end{pmatrix},\ens
and using the definition of $P_k(\ga)$ in \eqref{1001g}  these calculations can be summarized as
\bee \label{1101a} \beta_k(\ga)= \overline{P_k(\ga)} \, \alpha_k(\ga).\ene

Now, in terms of the vector $\alpha_k(\ga),$ the result in \eqref{1301b} means that
\bes   \sum_{\lambda\in \Lambda_{k+1}} | \la F, \ml \Phi_{k+1}\ra|^2    & = & \mu_{\hg}(V_{k+1}) \, \int_{V_{k}}
\alpha_k(\ga)^* \alpha_k(\ga)\, d\ga.\ens
Using the assumption \eqref{172e} and \eqref{1101a}, it follows that
\bes & \ & \sum_{\lambda\in \Lambda_{k+1}} | \la F, \ml \Phi_{k+1}\ra|^2  =
\mu_{\hg}(V_{k+1}) \, \frac{\mu_{\hg}(V_k)}{\mu_{\hg}(V_{k+1})}\, \int_{V_{k}}
\alpha_k(\ga)^* \overline{P_k(\ga)^* P_k(\ga) } \alpha_k(\ga)\, d\ga \\ & = &
\mu_{\hg}(V_k) \,  \int_{V_{k}}
\beta_k(\ga)^* \beta_k(\ga)\, d\ga  =  \mu_{\hg}(V_k) \,  \int_{V_{k}}\left( |b_k(\ga)|^2 + \| c_k(\ga)\|^2 \right)\, d\ga    \\ & = &
\mu_{\hg}(V_k) \,  \int_{V_{k}} \bigg(
\bigg| \sum_{\omega\in \Lambda_k^\bot} F(\omega + \ga) \overline{\Phi_{k}(\omega + \ga)} \bigg|^2
+ \sum_{m=1}^{\rho_k}\bigg| \sum_{\omega\in \Lambda_k^\bot} F(\omega + \ga) \overline{\Psi_{k}^{(m)}(\omega + \ga)} \bigg|^2\bigg)\, d\ga \\ & = &
\sum_{\lambda\in \Lambda_{k}} | \la F, \ml \Phi_k\ra|^2+
\sum_{m=1}^{\rho_k}\sum_{\lambda\in \Lambda_{k}} | \la F, \ml \Psi_k^{(m)}\ra|^2,\ens
where (\ref{1001n}) and (\ref{1001p}) are used in the final step.\ep

The next lemma is a consequence of the assumptions (i) and (ii) of the general setup.

\bl \label{901b} For any $F\in C_c(\hg)$ and  any $\epsilon>0,$ there is a $K\in I$ such that for $k\ge K,$ $k \in I,$
\bes
(1-\epsilon) \, \|F\|^2 \le \sulk | \la F, \ml \Phi_k\ra|^2 \le
(1+\epsilon) \, \|F\|^2.\ens \el

\bp Given $F\in C_c(\hg),$ put $S:= \supp\,F.$
For $k\in I, \, \omega \in \lhk,$ let
\bes
S_{k, \omega}:= \{\ga \in V_k \, \big| \, \omega + \ga \in S\}.\ens
Note that by \eqref{901f},
\bee \notag S & = & S \cap \hg  =  S \cap \bigg[ \bigcup_{\omega\in \lhk} (\omega + V_k)\bigg]
=   \bigcup_{\omega\in \lhk} \left[ S \cap (\omega + V_k)\right]
\\ & = &  \bigcup_{\omega\in \lhk} \{\omega + \ga \, \big| \ga\in V_k, \, \omega + \ga  \in S\}
=  \bigcup_{\omega\in \lhk} \{\omega + \ga \, \big| \ga\in S_{k, \omega}\}
\label{1001f} =  \bigcup_{\omega\in \lhk} (\omega + S_{k, \omega}). \hspace{.9cm}
\ene Since the decomposition in \eqref{901f} is disjoint, \eqref{1001f} is clearly a disjoint decomposition of $S$
(up to a set of measure zero for both decompositions).
By Lemma \ref{11m},
\bes \sulk |\la F, \ml \Phi_k \ra|^2 = \mu_{\hg}(V_k) \int_{V_k}
\bigg|  \sulhk F(\omega + \ga) \ \overline{\Phi_k(\omega + \ga)} \bigg|^2\, d\ga.\ens
Note that in the integral we only get contributions for the $\ga \in V_k$ for which
there is an $\omega^\prime \in \lhk$ such that $\omega^\prime + \ga \in S,$ i.e., we only get
contributions for $\ga \in \bigcup_{\omega^\prime \in \lhk} S_{k, \omega^\prime}.$ Thus
\bes \sulk |\la F, \ml \Phi_k \ra|^2 & = & \mu_{\hg}(V_k) \int_{\big[\bigcup_{\omega^\prime \in \lhk} S_{k, \omega^\prime}\big]}
\bigg|  \sulhk F(\omega + \ga) \ \overline{\Phi_k(\omega + \ga)} \bigg|^2\, d\ga. \ens

Now, given any $\epsilon>0,$ take $K\in I$ satisfying (i) and (ii) of the general setup. Then, for $k\in I$ with $k\ge K,$ the sets
$S_{k, \omega}, \, \omega \in \Lambda_k^\bot,$ are disjoint (up to a set of measure zero). Indeed, if
$\ga \in S_{k, \omega} \cap S_{k, \omega^\prime}$ for some $\omega, \omega^\prime\in \Lambda_k^\bot,$ then
$\ga \in (-\omega + S) \cap (-\omega^\prime +  S);$
by the assumption \eqref{901g} and the fact that $\Lambda_k^\bot \subset \Lambda_K^\bot,$ this implies that
$\mu_{\widehat{G}}( S_{k, \omega} \cap S_{k, \omega^\prime}) = 0$ if $\omega \neq \omega^\prime.$ We can therefore
continue our calculation, and obtain that for $k\ge K, k\in I,$
\bee \label{1101h}  \sulk |\la F, \ml \Phi_k \ra|^2
& = &   \mu_{\hg}(V_k) \sum_{\omega^\prime \in \lhk}   \int_{ S_{k, \omega^\prime} }
\bigg|  \sulhk F(\omega + \ga) \ \overline{\Phi_k(\omega + \ga)} \bigg|^2\, d\ga        .\ene  Note that for a fixed $\omega^\prime$ in the ``outer sum'',
we only get a nonzero contribution in the  ``inner sum" over  $\omega \in \Lambda_k^\bot$ for
the choice $\omega = \omega^\prime.$ In fact, given $\omega^\prime \in
\Lambda_k^\bot,$ for $\omega \neq \omega^\prime$ any $\ga \in S_{k, \omega^\prime}$
will be outside $S_{k, \omega},$ meaning that $\omega + \ga \notin S,$ i.e.,
$F(\omega + \ga)=0.$ Therefore \eqref{1101h} simplifies to
\bes \sulk |\la F, \ml \Phi_k \ra|^2
& = &   \mu_{\hg}(V_k) \sum_{\omega^\prime \in \lhk}   \int_{ S_{k, \omega^\prime} }
\big|  F(\omega^\prime + \ga) \Phi_k(\omega^\prime + \ga) \big|^2\, d\ga        \\ & = & \mu_{\hg}(V_k) \sum_{\omega^\prime \in \lhk}   \int_{ \omega^\prime + S_{k, \omega^\prime} }
\big|  F( \ga) \Phi_k(\ga) \big|^2\, d\ga  =
\mu_{\hg}(V_k)   \int_{ S }
\big|  F( \ga) \Phi_k(\ga) \big|^2\, d\ga,
\ens where the last step again used that the union in \eqref{1001f} is disjoint.
Our choice of $K$ and the assumption \eqref{901a} now implies that
\bes (1-\epsilon)\, \int_S |F(\ga)|^2\, d\ga \le  \sulk |\la F, \ml \Phi_k \ra|^2
\le (1+\epsilon)\, \int_S |F(\ga)|^2\, d\ga;\ens since
$\int_S |F(\ga)|^2\, d\ga= \| F\|^2,$ this completes the proof. \ep

We are now ready to state the main result, the unitary extension principle for LCA groups.

\bt \label{1001h} In addition to the assumptions \eqref{901g}--\eqref{1001a} in the general setup, assume that for $k\in I,$ the
matrix-valued function $P_k$ in \eqref{1001g} satisfies \eqref{172e}, i.e.,
\bes P_k(\ga)^* P_k(\ga)=d_k I_{d_k}, \, a.e. \, \ga \in V_k.\ens Then with
$\Psi_k^{(m)}, \, k\ge k_0, \, m=1, \dots, \rho_k,$ defined as in \eqref{1001b},
the collection of functions \bee \label{1001k} \{\ml \Phi_{k_0}\}_{\lambda\in \Lambda_{k_0}} \bigcup \,
\{\ml \Psi_k^{(m)}\}_{k\ge k_0, \lambda\in \Lambda_{k}, m=1, \dots, \rho_k}\ene forms a tight frame for $\lthg$ with frame bound $1$. \et

\bp Given any $\epsilon>0$ (the role hereof will be clear later), choose $K\in I$ such that \eqref{901g} and \eqref{901a} hold. Taking any $F\in C_c(\hg)$ and applying Lemma \ref{901bf} repeatedly, we obtain that for
$k\ge K, k\in I,$
\bes  \sum_{\lambda\in \Lambda_{k}} | \la F, \ml \Phi_{k}\ra|^2= \sum_{\lambda\in \Lambda_{k_0}} | \la F, \ml \Phi_{k_0}\ra|^2+\sum_{\ell= k_0}^{k-1}\sum_{m=1}^{\rho_\ell}
\sum_{\lambda\in \Lambda_{\ell}} | \la F, \ml \Psi_\ell^{(m)}\ra|^2.\ens
Using Lemma \ref{901b}, it follows that
\bee \label{3.21aaa} (1-\epsilon)\, \|F\|^2 \le \sum_{\lambda\in \Lambda_{k_0}} | \la F, \ml \Phi_{k_0}\ra|^2+\sum_{\ell= k_0}^{k-1}\sum_{m=1}^{\rho_\ell}
\sum_{\lambda\in \Lambda_{\ell}} | \la F, \ml \Psi_\ell^{(m)}\ra|^2\le (1+\epsilon)\, \|F\|^2.\ene
In the case where $I= \{k\}_{k= k_0}^\infty$, we can now let $k\to \infty,$ which yields
that \bes (1-\epsilon)\, \|F\|^2 \le \sum_{\lambda\in \Lambda_{k_0}} | \la F, \ml \Phi_{k_0}\ra|^2+\sum_{\ell\in I}\sum_{m=1}^{\rho_\ell}
\sum_{\lambda\in \Lambda_{\ell}} | \la F, \ml \Psi_\ell^{(m)}\ra|^2\le (1+\epsilon)\, \|F\|^2.\ens
Since $\epsilon>0$ is arbitrary, we conclude that
\bee \label{1101m}  \|F\|^2 = \sum_{\lambda\in \Lambda_{k_0}} | \la F, \ml \Phi_{k_0}\ra|^2+\sum_{\ell\in I}\sum_{m=1}^{\rho_\ell}
\sum_{\lambda\in \Lambda_{\ell}} | \la F, \ml \Psi_\ell^{(m)}\ra|^2.\ene
Due to the fact that $F$ is an arbitrary function in the dense subspace $C_c(\hg)$ of $\lthg,$ we infer that \eqref{1101m} holds for all $F\in \lthg.$

In the case where $I= \{k\}_{k= k_0}^{k_1},$ we simply take $k=k_1$ in \eqref{3.21aaa}. This eventually leads to the conclusion that the collection
$\{\ml \Phi_{k_0}\}_{\lambda\in \Lambda_{k_0}} \bigcup \,
\{\ml \Psi_k^{(m)}\}_{k = k_0, \dots, k_1 -1, \lambda\in \Lambda_{k}, m=1, \dots, \rho_k}$
is a tight frame for $\lthg$ with frame bound $1$.
\ep

Applying the inverse Fourier transform to the system in \eqref{1001k} immediately
leads to a generalized shift-invariant system forming a frame for $L^2(G)$:

\bc Under the assumptions in Theorem {\rm \ref{1001h}}, the generalized shift-invariant system
 \bes
 \{T_\lambda {\cal F}^{-1} \Phi_{k_0}\}_{\lambda\in \Lambda_{k_0}} \bigcup \,
\{T_\lambda {\cal F}^{-1} \Psi_k^{(m)}\}_{k\ge k_0, \lambda\in \Lambda_{k}, m=1, \dots, \rho_k}\ens forms a tight frame for $L^2(G)$ with frame bound $1$.\ec

The following example shows how the classical version of the UEP appears via Theorem \ref{1001h}. In particular, we will see that the condition of
\eqref{172e} holding for all values of $k$ reduces to just one matrix condition in this case. This is not surprising as Theorem \ref{1001h} on LCA groups is established for a much more general
nonstationary setting, whereas the classical UEP on the real line deals with stationary tight frames.

\bex \label{eg3.7} Consider an integer $a>1$ and the associated
{\em scaling operator} $D_a: \ltr \to \ltr, \, (D_af)(x):= a^{1/2}f(ax), \, x \in \mr.$ Given a function $\phi\in \ltr$
satisfying that
\begin{equation}
\label{3.30aaa}
 \lim_{\gamma \rightarrow 0} \widehat{\phi}(\gamma) = 1,
\end{equation}
let
$\Phi_k:= \widehat{D_a^k\phi}, \, k\in \mz.$ Also, let $\Lambda_k:=a^{-k}\mz.$ Then
$\Lambda_k^\bot=a^k \mz,$ with the fundamental domain $V_k=[0, a^k).$ Notice that
for this particular case the disjoint decomposition of $\Lambda_k^\bot$
in terms of $ \Lambda_{k+1}^\bot$, see   \eqref{1001d}, takes the form
$\Lambda_k^\bot= [0 + \Lambda_{k+1}^\bot]+ [a^k + \Lambda_{k+1}^\bot]+ \cdots
+ [(a-1)a^k+\Lambda_{k+1}^\bot].$ That is, $d_k=a,$ and
$\{\nu_{k, \ell}\}_{\ell=1}^{d_k}= \{ (\ell-1)\, a^k\}_{\ell =1}^{a}.$

In the general setup we consider the scaling relation
\bee \label{1201c} \Phi_{-1}(\ga)= H_0(\ga) \Phi_0(\ga), \, \ga \in \mr, \ene for a function $H_0\in L^\infty(V_0)= L^\infty([0,1)),$ extended to
a 1-periodic function.  The equation means that
${\cal F}D_a^{-1}\phi(\ga)=H_0(\ga) {\cal F}\phi(\ga)$ or
$\sqrt{a}\, \widehat{\phi}(a\ga)= H_0(\ga) \widehat{\phi}(\ga);$ this is the classical scaling equation in wavelet analysis, except that the factor $1/\sqrt{a}$ is usually absorbed in the function $H_0.$ It follows  that for any $k\in \mz,$
$\sqrt{a} \, \widehat{\phi}(a^k\ga)= H_0(a^{k-1}\ga) \widehat{\phi}(a^{k-1}\ga);$ translated back to the functions $\Phi_k$ this means that $\Phi_{-k}(\ga)= H_0(a^{k-1}\ga) \Phi_{-k+1}(\ga),$ or
\bes
\Phi_k(\ga)= H_0(a^{-k-1}\ga) \Phi_{k+1}(\ga), \, \ga \in \mr.\ens
That is, the single scaling equation \eqref{1201c} implies that we have the refinement equation $\Phi_k(\ga)= H_{k+1}(\ga) \Phi_{k+1}(\ga)$ for all levels, with
$H_{k+1}(\ga)= H_0(a^{-k-1}\ga).$
As for the remaining assumptions (\ref{901g}) and (\ref{901a}) in the general setup, (\ref{901g}) follows from the fact that $\Lambda_k^\bot = a^k \mz$,
while (\ref{901a}) is a consequence of the condition (\ref{3.30aaa}) and the calculation
$$
 \mu_{\mr}(V_k)| \Phi_k(\ga)|^2 = a^k |a^{-k/2}\widehat{\phi}(a^{-k}\ga)|^2 = |\widehat{\phi}(a^{-k} \ga)|^2.
$$

Assuming now that we have chosen functions
$G_0^{(m)}\in L^\infty([0,1)),$ $m=1,\dots,\rho,$ for some integer $\rho \geq a -1,$ define the functions $G_{k+1}^{(m)}$ by
$G_{k+1}^{(m)}(\ga):= G_0^{(m)}(a^{-k-1}\ga), \, \ga \in \mr,$ for $m=1, \dots, \rho,$ i.e., we take $\rho_k =\rho$ for all $k \in \mz.$
Considering the entries in the first row of the $(\rho + 1) \times a$ matrix $P_k(\ga)$ in (\ref{1001g}), we see that for $\ga \in V_k=[0, a^k),$
\bes H_{k+1}(\ga + \nu_{k, \ell})= H_0(a^{-k-1}(\ga+ (\ell-1)\, a^{k}))
  = H_0(a^{-k-1}\ga+ \nu_{-1, \ell}). \ens
A similar calculation works for $G_k^{(m)},$ $m =1,\dots, \rho,$ so we conclude that for $\ga \in V_k,$
$P_k(\ga)=P_{-1}(a^{-k-1}\ga).$ Since $a^{-k-1}\ga$ runs through $V_{-1}=[0, a^{-1})$ when
$\ga$ runs through $V_k,$ we conclude that for all $k\in \mz$ the matrix equation
\eqref{172e} is equivalent to
$P_{-1}(\ga)^* P_{-1}(\ga)= a\, I_{a}, \, \mbox{a.e.} \, \ga \in V_{-1},$ which corresponds to the condition in the classical UEP on $\mr.$
Note that for $\lambda \in \Lambda_k= a^{-k}\mz,$ i.e., $\lambda= a^{-k}\ell$ for some $\ell \in \mz,$
\bes \ml \Phi_k(\ga)= e^{-2 \pi i a^{-k}\ell \ga} {\cal F} D_a^k\phi(\ga)=
{\cal F}T_{\ell a^{-k}}D_a^k\phi(\ga),\ens with a similar calculation yielding
$ \ml \Psi_k^{(m)}(\ga)= {\cal F}T_{\ell a^{-k}}D_a^k\psi^{(m)}(\ga), \, m=1\ldots, \rho,$ where
$\widehat{\psi^{(m)}}:= \Psi_0^{(m)}$ and $\Psi_k^{(m)} = \widehat{D_a^k \psi^{(m)}}.$
If the assumptions in Theorem \ref{1001h} are satisfied, this implies that
for any chosen $k_0\in \mz$ (with $I= \{k\}_{k=k_0}^\infty$), the collection of functions
\begin{eqnarray*}  & & \{\ml \Phi_{k_0}\}_{\lambda\in \Lambda_{k_0}} \bigcup \,
\{\ml \Psi_k^{(m)}\}_{k\ge k_0, \lambda\in \Lambda_{k}, m=1,\dots,\rho} \nonumber \\ & = &
\{ {\cal F}T_{\ell a^{-{k_0}}}D_a^{k_0}\phi\}_{\ell\in \mz} \bigcup \, \{ {\cal F}T_{\ell a^{-k}}D_a^k\psi^{(m)}\}_{k\ge k_0, \ell\in \mz,
m=1,\dots,\rho}\end{eqnarray*}
forms a tight frame for $\ltr$ with frame bound 1. Thus the system
\bes & &
\{ T_{\ell a^{-{k_0}}}D_a^{k_0}\phi\}_{\ell\in \mz} \bigcup \, \{ T_{\ell a^{-k}}D_a^k\psi^{(m)}\}_{k\ge k_0, \ell\in \mz, m=1,\dots,\rho} \\ & = &
\{ a^{k_0/2}\phi(a^{k_0}\cdot-\ell)\}_{\ell\in \mz} \bigcup \, \{a^{k/2}\psi^{(m)}(a^{k}\cdot-\ell) \}_{k\ge k_0, \ell\in \mz,
m=1,\dots,\rho}
\ens
forms a tight frame for $\ltr$ with frame bound 1.
\ep \enx

A more recent variant of the UEP, the {\it oblique extension principle} (OEP), was announced in 2001,
independently by two groups of researchers, namely, Chui, He and St\"ockler \cite{CHS}, as well
as Daubechies, Han, Ron and Shen \cite{DRoSh3}. The OEP is essentially equivalent to the UEP, but gives easier access to attractive constructions, e.g., constructions with a high number of vanishing moments. Theorem \ref{1001h} can be extended to an OEP version along the lines of the proofs in
\cite{CHS, DRoSh3}.

Without extra mathematical difficulties (except the ones that come from a  more cumbersome notation), Theorem \ref{1001h}
can be generalized to the multi-generator case. That is, we can replace each of the functions $\Phi_k, \, k\in I,$ in the general setup by
a vector of functions ${\bf \Phi}_k(\ga):= (\Phi_k^{(1)}(\ga), \dots, \Phi_k^{(r_k)} (\ga))^T$
and the condition $\Phi_k (\ga)= H_{k+1}(\ga) \Phi_{k+1}(\ga)$
by a matrix condition
$ {\bf \Phi}_k(\ga)= {\bf H}_{k+1}(\ga) {\bf \Phi}_{k+1}(\ga),$
where  ${\bf H}_{k+1}(\ga)$ is an $r_{k} \times r_{k+1}$ matrix.

Similarly, assume that for $k\in I$ we are given a $\rho_{k} \times r_{k+1}$
matrix-valued function ${\bf G}_{k+1}$ with entries in $ L^\infty(V_{k+1}),$ define the
functions
$ {\bf \Psi}_k(\ga):= (\Psi_k^{(1)}(\ga), \dots,\Psi_k^{(\rho_k)} (\ga))^T := {\bf G}_{k+1}(\ga) {\bf \Phi}_{k+1}(\ga).$
The
corresponding version of the matrix $P_k(\ga)$ in
\eqref{1001g} is
\bes {\bf P}_k(\ga):= \begin{pmatrix} {\bf H}_{k+1} (\ga + \nu_{k,1}) & \cdot & \cdot & \cdot & {\bf H}_{k+1} (\ga + \nu_{k,d_k}) \\
{\bf G}_{k+1} (\ga + \nu_{k, 1}) & \cdot & \cdot & \cdot & {\bf G}_{k+1} (\ga + \nu_{k,d_k})  \end{pmatrix},\ens which is now an
$(r_k + \rho_k) \times r_{k+1}d_k$ matrix. Finally, while keeping the condition (\ref{901g}), replace the technical condition \eqref{901a} by
the assumption that for every compact set $S$ in $\hg$ and any $\epsilon >0,$ there exists $K_2\in I$ such that for all $k\ge K_2, \, k\in I,$
\bes
\bigg| \mu_{\hg}(V_k) \, \sum_{m=1}^{r_k} |\Phi_k^{(m)}(\ga) |^2 -1 \bigg| \le \epsilon, \ \forall \ga\in S.\ens

Such a multi-generator setting is considered in \cite{GT2} for the UEP on $\mt^s$. Following the arguments in \cite{GT2}, the proof of
Theorem \ref{1001h} can be easily adapted to give a UEP on LCA groups for the multi-generator setting.

The multi-generator version of Theorem \ref{1001h} now says that if for every
$k\in I,$
$${\bf P}_k(\ga)^* {\bf P}_k(\ga)= d_k I_{r_{k+1}d_k}, \, \mbox{a.e.} \, \ga\in V_k,$$ then the collection of functions
\bes
\{\ml \Phi_{k_0}^{(m)}\}_{\lambda\in \Lambda_{k_0}, m=1, \dots, r_{k_0}} \bigcup \,
\{\ml \Psi_k^{(m)}\}_{k\ge k_0, \lambda\in \Lambda_{k}, m=1, \dots, \rho_k}\ens forms a tight frame for $\lthg$ with frame bound 1.

\section{B-spline generated systems}     \label{1901d}
In this section we will consider certain explicitly given functions on LCA groups
and verify some of the UEP conditions without
any restriction on the underlying LCA group. The remaining conditions will be verified directly in concrete cases.

The considered systems will be based on a version of the B-splines on LCA groups, as
introduced independently by Dahlke \cite{Dah} and Tikhomirov \cite{T} in 1994.
The results generalize the well-known construction by Ron and Shen \cite{RoSh2} of a tight frame based on B-splines on $G=\mr.$

Consider a nested sequence of lattices in the LCA group $G,$ as in \eqref{901k}.
We will assume that for each $k\in I,$
\bee \label{41229a}|\Lambda_{k+1}/\Lambda_k|=2\ene and apply Lemma \ref{801b} on each inclusion $\Lambda_k \subset \Lambda_{k+1}.$
In particular,  there is a $\nu_k\in \Lambda_k^\bot \setminus \Lambda_{k+1}^\bot$ such that we have the disjoint splitting
\bee \label{1901h}  \Lambda_k^\bot = \Lambda_{k+1}^\bot \cup (\nu_k + \Lambda_{k+1}^\bot).\ene
Similarly, there exists $\eta_k\in \Lambda_{k+1} \setminus \Lambda_k$ such that we have the disjoint splitting
\bee \label{1901j} \Lambda_{k+1} = \Lambda_k \cup (\eta_k + \Lambda_k).\ene
In the entire section we let $V_k$ denote an arbitrary fundamental domain associated with the lattice $\lak^\bot.$ Now,
fix $k^\prime\ge k_0$ and let $Q_{k^\prime}$ denote a fundamental domain associated with the lattice $\Lambda_{k^\prime}$ in $G.$ It follows from Lemma \ref{801b} that we
successively can construct fundamental domains  $Q_k$ for the lattices $\Lambda_k$ such that for all $k\le k^\prime$  we have the disjoint splitting
\bee \label{1901n}  Q_k = Q_{k+1} \cup (\eta_k + Q_{k+1}).\ene
We will assume that \eqref{1901n} holds for all $k\in I$ (we will later verify
this condition explicitly in concrete cases). On the other hand,
we note that it is actually possible to construct ``badly chosen" fundamental domains $Q_k$ such that no fundamental domain $Q_{k+1}$ will satisfy
\eqref{1901n} with a disjoint splitting, regardless of the choice of $\eta_k
\in \Lambda_{k+1} \setminus \Lambda_k.$ We will verify \eqref{1901n} directly for the concrete cases to be considered later.
\bl \label{1901r} Assume that \eqref{41229a} holds and choose  $\nu_k, \eta_k$
as in \eqref{1901h} and \eqref{1901j}. Then
$(\eta_k, \nu_k)=-1$ for all $k\in I.$ \el

\bp Since $\nu_k\in \Lambda_k^\bot \setminus \Lambda_{k+1}^\bot,$ we see that
$\nu_k+\nu_k \in \Lambda_k^\bot \setminus (\nu_k+ \Lambda_{k+1}^\bot);$ using the decomposition \eqref{1901h}, we conclude that $\nu_k+\nu_k\in \Lambda_{k+1}^\bot.$
Since $\eta_k\in \Lambda_{k+1},$ it follows that
$1= (\eta_k, \nu_k+ \nu_k)= (\eta_k, \nu_k)^2,$ i.e., $(\eta_k, \nu_k)= \pm 1.$ But if
$(\eta_k, \nu_k)=1,$ then for any $\lambda\in \Lambda_k, \, (\eta_k+\lambda, \nu_k)=1$ and $(\lambda, \nu_k)=1;$ by the decomposition \eqref{1901j}, this shows that
$(\lambda^\prime, \nu_k)=1$ for all $\lambda^\prime \in \Lambda_{k+1},$ i.e.,
$\nu_k \in \Lambda_{k+1}^\bot.$ This contradicts that $\nu_k\in \Lambda_k^\bot \setminus \Lambda_{k+1}^\bot,$ and we conclude that $(\eta_k, \nu_k)=-1.$\ep

Our constructions will be based on the following definition of B-splines on LCA groups.

\bd Consider a sequence of nested lattices $\{\Lambda_k\}_{k=k_0}^\infty$ in $G,$ with associated fundamental domains $Q_k, k\ge k_0.$ For $N\in \mn,$ define the B-spline of $N$th order
at level $k, k\ge k_0,$ by the $N$-fold convolution
\bes
\phi_{k,N}(x):= \mu_G(Q_k)^{-N+1/2} \chi_{Q_k} * \cdots * \chi_{Q_k}(x), \, x\in G.\ens \ed

We have included the factor $\mu_G(Q_k)^{-N+1/2}$ in the definition of the B-spline in order to avoid a later renormalization. Note also that the index $k$ refers
to the level of the  sets $Q_k$ within the scale of fundamental domains associated with the lattices $\Lambda_k.$  We will consider a fixed choice of $N\in \mn$ and consider the functions $\Phi_{k,N}, \,
k\ge k_0,$ defined by
\bee \label{1901t} \Phi_{k,N}(\ga):= \widehat{\phi_{k,N}}(\ga)= \mu_G(Q_k)^{-N+1/2}
\left( \int_{Q_k} (-x, \ga)\, dx\right)^N, \, \ga \in \hg,\ene
where the Fourier transform formula (\ref{2.1aaa}) is used to obtain the final expression.

We now show that for any given $N\in \mn,$ the functions $\Phi_{k,N}$ satisfy the refinement equation \eqref{1001a}; we suppress the dependence of $N$ in the notation for the associated filters $H_{k+1}.$ For the one-dimensional elementary groups this result was also proved in \cite{S94}.

\bl \label{2001e} Assume that  \eqref{41229a} and \eqref{1901n} are satisfied. Fix $N\in \mn.$ The functions $\Phi_{k,N},$ $k\in I,$ defined by \eqref{1901t} satisfy the refinement equation
\bee \label{1901u} \Phi_{k,N}(\ga) = H_{k+1}(\ga) \Phi_{k+1,N}(\ga), \, \ga\in \hg,\ene where $H_{k+1}\in L^\infty(V_{k+1})$ is given by
\bee \label{1901uu} H_{k+1}(\ga)= \frac1{2^{N-1/2}}\, (1+ (-\eta_k, \ga))^N, \, \ga \in \hg.\ene \el
\bp Using the disjoint splitting in \eqref{1901n} and a change of variable,
\bes \Phi_{k,N}(\ga) & = &
\mu_G(Q_k)^{-N+1/2}
\left( \int_{Q_{k+1}} (-x, \ga)\, dx+  \int_{\eta_k +Q_{k+1}} (-x, \ga)\, dx \right )^N \\ & = &\mu_G(Q_k)^{-N+1/2}
 (1+ (-\eta_k, \ga))^N \left( \int_{Q_{k+1}} (-x, \ga)\, dx \right)^N \\ & = & \left(\frac{\mu_G(Q_k)}{\mu_G(Q_{k+1})}\right)^{-N+1/2}
 (1+ (-\eta_k, \ga))^N \mu_G(Q_{k+1})^{-N+1/2}\left( \int_{Q_{k+1}} (-x, \ga)\, dx \right)^N. \ens
Note that by Lemma \ref{801b}(i),
$2= | \Lambda_{k+1}/\Lambda_k| = \frac{s(\Lambda_k)}{s(\Lambda_{k+1})}
= \frac{\mu_G(Q_k)}{\mu_G(Q_{k+1})};$ thus, we see that \eqref{1901u} is satisfied with
the function $H_{k+1}$ defined as in \eqref{1901uu}. Clearly $H_{k+1}$ is bounded; that
$H_{k+1}$ is periodic follows from the fact that for $\omega \in \Lambda_{k+1}^\bot,
\ga \in \hg,$
\bqs H_{k+1}(\omega + \ga)= \frac{1}{2^{N-1/2}}
 (1+ (-\eta_k, \omega +\ga))^N =\frac{1}{2^{N-1/2}}
 (1+ (-\eta_k, \ga))^N =H_{k+1}(\ga). \eqs

We will now provide a condition under which \eqref{901a} holds.

\bl \label{2001a} Assume that \eqref{1901n} holds and fix $N\in \mn.$ Then the following hold:
\bei \item[{\rm (i)}] Let $\delta\in (0,1)$ be given. Assume that for some $k\in I$ and
some $\ga \in \hg,$
\bee \label{2001b} | (-x, \ga)-1| \le \delta, \, \forall x\in Q_k.\ene Then
\bee \label{2001c}  | \mu_{\hg}(V_k) \, |\Phi_{k,N}(\ga) |^2 -1 | \le 1 - (1-\delta)^{2N}. \ene
\item[{\rm (ii)}] Assume that for every $\delta \in (0,1)$ and any compact set $S$ in $\hg,$ there exists $k\in I$ such that the inequality \eqref{2001b} holds for all
    $\ga \in S$ and $x\in Q_k.$ Then the condition \eqref{901a} is satisfied.
\eni
\el

\bp First, using \eqref{1901t} and that $\mu_{\hg}(V_k) \mu_G (Q_k)=1,$
\bes \mu_{\hg}(V_k) |\Phi_{k,N}(\ga)|^2 = \mu_{\hg}(V_k) \mu_G(Q_k)^{-2N+1}
\bigg| \int_{Q_k} (-x, \ga)\, dx\bigg|^{2N} =
\bigg| \frac1{\mu_G(Q_k)} \int_{Q_k} (-x, \ga)\, dx\bigg|^{2N}. \ens
It follows that
\bes \mu_{\hg}(V_k) |\Phi_{k,N}(\ga)|^2 \le \bigg( \frac1{\mu_G(Q_k)} \int_{Q_k} |(-x, \ga)|\, dx\bigg)^{2N}=1,\ens and therefore
\bee \label{2001d} | \mu_{\hg}(V_k) |\Phi_{k,N}(\ga) |^2 -1 | =
 1-\mu_{\hg}(V_k) |\Phi_{k,N}(\ga) |^2= 1- \bigg| \frac1{\mu_G(Q_k)} \int_{Q_k} (-x, \ga)\, dx\bigg|^{2N}.\ene Using the triangle inequality, the assumption \eqref{2001b} implies that
 \bes \bigg| \,  \bigg| \frac1{\mu_G(Q_k)} \int_{Q_k} (-x, \ga)\, dx \bigg|  -1 \bigg|
& \le & \bigg|    \frac1{\mu_G(Q_k)} \int_{Q_k} (-x, \ga)\, dx -1 \bigg| \\
& = & \bigg|    \frac1{\mu_G(Q_k)} \int_{Q_k} [(-x, \ga)-1]\, dx  \bigg| \\
& \le & \frac1{\mu_G(Q_k)} \int_{Q_k} |(-x, \ga)-1|\, dx \le \delta,\ens and therefore
$ \big|  \frac1{\mu_G(Q_k)} \int_{Q_k} (-x, \ga)\, dx  \big| \ge 1-\delta.$
Via  \eqref{2001d} we conclude that
$ | \mu_{\hg}(V_k) |\Phi_{k,N}(\ga) |^2 -1 | \le 1 - (1-\delta)^{2N},$ which proves (i).

In order to prove (ii), given $\epsilon>0,$ choose $\delta \in (0,1)$ such that
$1-(1-\delta)^{2N} \le \epsilon.$ Then by (\ref{2001c}), $| \mu_{\hg}(V_k) |\Phi_{k,N}(\ga) |^2 -1 |\le \epsilon$ for all $\ga \in S.$ Since $Q_{k+1} \subset Q_k$ for all $k\in I,$
the same holds with $k$ replaced by any $k^\prime$ for which $k^\prime \ge k.$ \ep

Now, fix $N\in \mn$ and consider the functions $\Phi_{k,N}, \, k \in I.$
Take the functions $H_{k´+1}, \, k\in I,$ as in \eqref{1901uu}; our task is to find functions
$G_{k+1}^{(1)}, \dots G_{k+1}^{(\rho_k)}\in L^\infty(V_{k+1})$ such that the matrices
\bee \label{1001ga} P_k(\ga):= \begin{pmatrix} H_{k+1} (\ga ) &  H_{k+1} (\ga + \nu_{k}) \\
G_{k+1}^{(1)} (\ga ) &  G_{k+1}^{(1)} (\ga + \nu_{k}) \\ \cdot & \cdot  \\
 \cdot & \cdot  \\
G_{k+1}^{(\rho_k)} (\ga ) & G_{k+1}^{(\rho_k)} (\ga + \nu_{k})
\end{pmatrix}, \ \ga \in V_k,\ene satisfy the matrix condition in the UEP, i.e.,
\bee \label{1001gb} P_k(\ga)^*P_k(\ga)= 2 I_2, \, \mbox{a.e.} \, \ga \in V_k.\ene
As we will see, the LCA case is technically more involved than the case of the real line.
We will consider the cases $N=1$  and $N=2M, \, M\in \mn.$
In both instances, we take $\rho_k$ to be $N,$ and for the case $N=1,$ we write $G_{k+1}^{(1)}$ simply as $G_{k+1}.$

\bpr \label{41230a} Assume that  \eqref{41229a} and \eqref{1901n} are satisfied.
For every $k \geq k_0,$ consider the following two choices of $\Phi_k$ and associated filters:
\bei
\item[{\rm (i)}] Let $N=1$ and consider
$\Phi_{k}:=\widehat{\phi_{k,1}},$ with the associated filter  $H_{k+1}$
in  \eqref{1901uu}.  Define the function $G_{k+1} \in L^{\infty}(V_{k+1})$ by
\bee \label{42101e} G_{k+1}(\ga):=   \frac1{\sqrt{2}}\, (1- (-\eta_k, \ga)), \, \gamma \in \hg. \ene
\item[{\rm (ii)}] Given $M\in \mn,$ let
$\Phi_k:=\widehat{\phi_{k,2M}},$ with  the associated filter  $H_{k+1}$
in  \eqref{1901uu}.  Define the functions $G_{k+1}^{(1)}, \dots G_{k+1}^{(2M)}\in L^\infty(V_{k+1})$ by
\bee \label{2101e} G_{k+1}^{(m)}(\ga):= \frac1{2^{2M-1/2}}\, \sqrt{\begin{pmatrix}  2M \\ m \end{pmatrix}}(1+ (-\eta_k, \ga))^{2M-m} (1- (-\eta_k, \ga))^{m},
\, \ga \in \hg.\ene
\eni
Then, in both cases  the matrix $P_k(\ga)$ in \eqref{1001ga} satisfies the UEP condition \eqref{1001gb}.
\epr

\bp For the proof of (i), by Lemma \ref{2001e},
the function $\Phi_{k}$ satisfies the refinement equation (\ref{1901u}) with
$H_{k+1}(\ga)= \frac1{\sqrt{2}}\, (1+ (-\eta_k, \ga)), \, \ga \in \hg.$
Now, take $G_{k+1}$ as in \eqref{42101e}. Then, by direct calculation,
\bes & & | H_{k+1}(\ga)|^2 + | G_{k+1}(\ga)|^2 \\ & = & \frac12 \Big((1+(-\eta_k, \ga))\overline{(1+(-\eta_k, \ga))} + (1-(-\eta_k, \ga))\overline{(1-(-\eta_k, \ga))}\Big)=2. \ens
Similarly, using that $(\eta_k, \nu_k)=-1$ by Lemma \ref{1901r},  we arrive at
\bes & \ & H_{k+1}(\ga) \overline{H_{k+1}(\ga + \nu_k)}+ G_{k+1}(\ga) \overline{G_{k+1}(\ga + \nu_k)} \\ & = & \frac12 \big( (1+ (-\eta_k, \ga)) (1+ (\eta_k, \ga+ \nu_k)) + (1- (-\eta_k, \ga)) (1- (\eta_k, \ga+ \nu_k)) \big)
  \\ & = & \frac12 \big( (1+ (-\eta_k, \ga)) (1- (\eta_k, \ga)) + (1- (-\eta_k, \ga)) (1+ (\eta_k, \ga))    \big) =0.
\ens
For the proof of (ii), clearly $G_{k+1}^{(m)}\in L^\infty(V_{k+1})$ for $m=1,\dots, 2M.$
We first show that
\bee \label{2101c} | H_{k+1}(\ga)|^2 + \sum_{m=1}^{2M} |G_{k+1}^{(m)}(\ga)|^2=2, \, \ga \in \hg.\ene
We will now use  the elementary identity
$ |1+z|^2 + |1-z|^2=4,$ which is valid for all $z\in \mt.$
Taking $z=(-\eta_k, \ga)$
and raising both sides  to the power of $2M$ yields that
 \bes \frac1{2^{4M-1}} \left(|1+(-\eta_k, \ga)|^2 + |1-(-\eta_k, \ga)|^2 \right)^{2M}=2,\ens
 which by the binomial theorem means that
\bes \frac1{2^{4M-1}}\sum_{m=0}^{2M}\begin{pmatrix}  2M \\ m \end{pmatrix}
|1+(-\eta_k, \ga)|^{4M-2m}  |1-(-\eta_k, \ga)|^{2m}=2.\ens This immediately yields that
\eqref{2101c} is satisfied; in fact, using \eqref{1901uu} and \eqref{2101e},
\bes & \ &| H_{k+1}(\ga)|^2 + \sum_{m=1}^{2M} |G_{k+1}^{(m)}(\ga)|^2 = \\ & \ &
\frac1{2^{4M-1}}\, |(1+ (-\eta_k, \ga))|^{4M}+ \frac1{2^{4M-1}}\sum_{m=1}^{2M}
\begin{pmatrix}  2M \\ m \end{pmatrix}
|(1+ (-\eta_k, \ga))|^{4M-2m} |(1- (-\eta_k, \ga))|^{2m} = 2.
\ens

We now have to show that
\bee \label{2101ca}  H_{k+1}(\ga) \overline{H_{k+1}(\ga+\nu_k)} + \sum_{m=1}^{2M} G_{k+1}^{(m)}(\ga) \overline{G_{k+1}^{(m)}(\ga+\nu_k)}=0, \, \ga \in \hg.\ene
Use the identity
$(1+z)(1-\overline{z}) + (1-z)(1+\overline{z})=0, z\in \mt,$
with $z=(-\eta_k, \ga);$ then
\bes \left((1+(-\eta_k, \ga))(1-\overline{(-\eta_k, \ga)}) + (1-(-\eta_k, \ga))(1+\overline{(-\eta_k, \ga)})\right)^{2M}=0,\ens which yields that
\bee \label{2101f} \sum_{m=0}^{2M}
\begin{pmatrix}  2M \\ m \end{pmatrix}
(1+ (-\eta_k, \ga))^{2M-m} (1- \overline{(-\eta_k, \ga)})^{2M-m}(1- (-\eta_k, \ga))^{m} (1+ \overline{(-\eta_k, \ga)})^{m}=0. \, \hspace{.3cm} \ \ene Inserting the expressions for $H_{k+1}$ and
$G_{k+1}^{(m)}$ into the left-hand side of \eqref{2101ca} yields
\bes  & & H_{k+1}(\ga) \overline{H_{k+1}(\ga+\nu_k)} + \sum_{m=1}^{2M} G_{k+1}^{(m)}
\overline{G_{k+1}^{(m)}(\ga+\nu_k)}  \\ & = & \frac1{2^{4M-1}}\, (1+ (-\eta_k, \ga))^{2M}(1+ \overline{(-\eta_k, \ga+\nu_k)})^{2M} \\
& \ & \  + \ \frac1{2^{4M-1}} \sum_{m=1}^{2M}
\begin{pmatrix}  2M \\ m \end{pmatrix}
(1+ (-\eta_k, \ga))^{2M-m} (1- (-\eta_k, \ga))^{m}
\\
& \ & \hspace{3.6cm}
\times
(1+ \overline{(-\eta_k, \ga + \nu_k)})^{2M-m} (1- \overline{(-\eta_k, \ga+ \nu_k)})^{m}.\ens
Using again that by Lemma \ref{1901r} $(\eta_k, \nu_k)=-1,$ it follows immediately from
\eqref{2101f} that this expression equals zero, as desired.
\ep

Applying the UEP to the setting in Proposition \ref{41230a}(ii), we obtain the following generalization to LCA groups of the B-spline tight frames on $\mr$
by Ron and Shen \cite{RoSh2}.

\bpr \label{41229c} Assume that  \eqref{901g},  \eqref{41229a} and \eqref{1901n} are satisfied, and that for each $\delta \in (0,1)$ and each compact set $S\subset \hg,$
there exists $k\in I$ such that the inequality \eqref{2001b} holds for all
    $\ga \in S$ and $x\in Q_k.$ Given $M\in \mn,$ for every $k \geq k_0,$ let
$\Phi_k:=\widehat{\phi_{k,2M}},$ with  the associated filter  $H_{k+1}$
in  \eqref{1901uu}.  Define the functions $G_{k+1}^{(m)}, m=1, \dots, 2M,$  by
\eqref{2101e}. Then the collection of functions
\bes
\{\ml \Phi_{k_0}\}_{\lambda\in \Lambda_{k_0}} \bigcup \,
\{\ml \Psi_k^{(m)}\}_{k\ge k_0, \lambda\in \Lambda_{k}, m=1, \dots, 2M}\ens forms a tight frame for $\lthg$ with frame bound $1$;
equivalently,
\bes
\{T_\lambda  \phi_{k_0,2M}\}_{\lambda\in \Lambda_{k_0}} \bigcup \,
\{T_\lambda  \psi_k^{(m)}\}_{k\ge k_0, \lambda\in \Lambda_{k}, m=1, \dots, 2M}\ens forms a tight frame for $L^2(G)$ with frame bound $1$. \epr

For the first-order B-splines in Proposition \ref{41230a}(i), we can of course formulate a UEP-result of exactly
the same type as in Proposition \ref{41229c};
in this particular case we actually obtain  an orthonormal basis. In order to show this, it is sufficient to establish that
$\| \Phi_k \|=\| \Psi_k\|=1$ for all $k \in I,$ so let us do that. By the definition of $\Phi_k$ and Parseval's identity, we immediately see that $\| \Phi_k\|=1$ for all $k\in I.$ Now, by
the definition of $\Psi_k$ in (\ref{1001b}), using (\ref{1901t}), (\ref{42101e}) and then the disjoint splitting (\ref{1901n}), for $\ga \in \hg,$
\bes \Psi_k(\ga)  & = &
\frac1{\sqrt{2}}\, (1- (-\eta_k, \ga)) \frac1{\sqrt{\mu_G (Q_{k+1})}}\int_{Q_{k+1}} (-x, \ga)\, dx  \\ & = & \frac1{\sqrt{2}}  \frac1{\sqrt{\mu_G (Q_{k+1})}}
\bigg(  \int_{Q_{k+1}} (-x, \ga)\, dx -\int_{Q_{k+1}} (-x-\eta_k, \ga)\, dx    \bigg)   \\ & = & \frac1{\sqrt{2}} \frac1{\sqrt{\mu_G (Q_{k+1})}}\bigg( \int_{Q_{k+1}} (-x, \ga)\, dx - \int_{\eta_k+Q_{k+1}} (-x, \ga)\, dx \bigg) \\
 & = &  \frac1{\sqrt{2}} \frac1{\sqrt{\mu_G (Q_{k+1})}}\left(   \widehat{\chi_{Q_{k+1}}}(\ga) - \widehat{\chi_{\eta_k + Q_{k+1}}}(\ga)  \right). \ens
Thus,
\bes {\cal F}^{-1}\Psi_k(x)=\frac1{\sqrt{2}}\frac1{\sqrt{\mu_G (Q_{k+1})}}\left(  \chi_{Q_{k+1}}(x) - \chi_{\eta_k + Q_{k+1}}(x)  \right), \, x\in G. \ens Using again the disjoint splitting in \eqref{1901n},
it now follows that
\bes \| \Psi_k\|^2 = \| {\cal F}^{-1}\Psi_k\|^2 = \frac12 \frac1{\mu_G(Q_{k+1})} \left( \mu_G(Q_{k+1})+ \mu_G(\eta_k+Q_{k+1})  \right) =1, \ens
as claimed.

We will now provide completely explicit constructions of frames  for a number of LCA groups. We will show that all the technical assumptions can be fulfilled in these cases.
Then the filters $G_{k+1}^{(m)}$ given in Proposition \ref{41230a} would lead to the desired constructions of $\Psi_k^{(m)}$ via Proposition \ref{41229c} and its above variant.

\bex \label{41228a} Consider the group $G= \mz,$ with dual group $\hg= \mt.$ Given
$M\in \mn,$ we will consider the lattices $ \label{4483} \Lambda_k:=2^{M-k} \mz, \, k=0, \dots, M;$ as associated fundamental domains, we take
$\label{4485} Q_k= \{0, \dots, 2^{M-k}-1\}.$
Note that $\Lambda_0=2^M\mz$ and $\Lambda_M=\mz=G.$ Clearly,
for $k=0, \dots, M-1,$
$\Lambda_k\subset \Lambda_{k+1}$
and $| \Lambda_{k+1}/\Lambda_k|=2;$ furthermore,
$\label{4484} \Lambda_{k+1}= \lak \cup \left(2^{M-k-1} + \lak\right)$, i.e.,
\eqref{1901j} holds with $\eta_k=2^{M-k-1}.$ Note that
\bes
Q_k =Q_{k+1} \cup \left(2^{M-k-1} + Q_{k+1}\right), \ens
i.e., \eqref{1901n} holds for all $k=0, \dots, M-1.$
Also, for $k=0, \dots,M,$
$\label{4484a} \lak^\bot=\frac1{2^{M-k}}\, \mz_{2^{M-k}};$ the set
$ \label{4484b}  V_k=[0, \frac1{2^{M-k}})$ is a fundamental domain associated with
the lattice $\lak^\bot.$ Observe that for $k =0, \dots, M-1,$
$$
\lak^\bot= \lako^\bot + \big( \frac1{2^{M-k}} + \lako^\bot\big),
$$
i.e.,  \eqref{1901h} holds with $\nu_k= \frac1{2^{M-k}}\in \lak^\bot \setminus \lako^\bot.$

We now verify the conditions \eqref{2001b} and \eqref{901g}. To this end, let $0<\delta <1$ and consider a compact set $S$ in $\hg=\mt.$ Then, taking $k=M$ and noticing that
$Q_M=\{0\},$ the condition \eqref{2001b} trivially holds. Also, $\Lambda_M^\bot= \{0\},$
so \eqref{901g} clearly holds.

\begin{figure}
\centering
\begin{tabular}{c}
\resizebox{2.5in}{2.2in}{\includegraphics{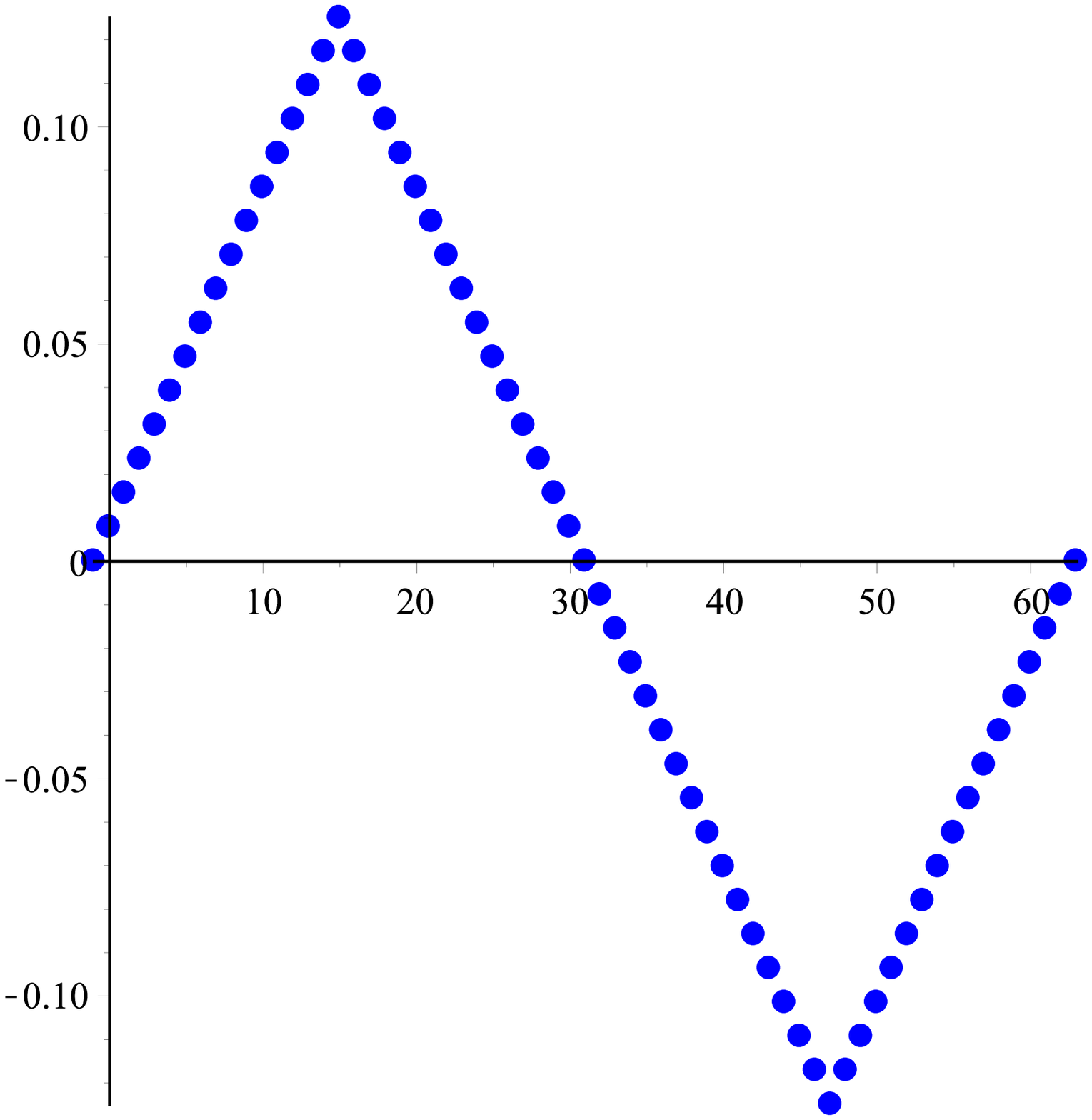}} \quad\quad
\resizebox{2.5in}{2.2in}{\includegraphics{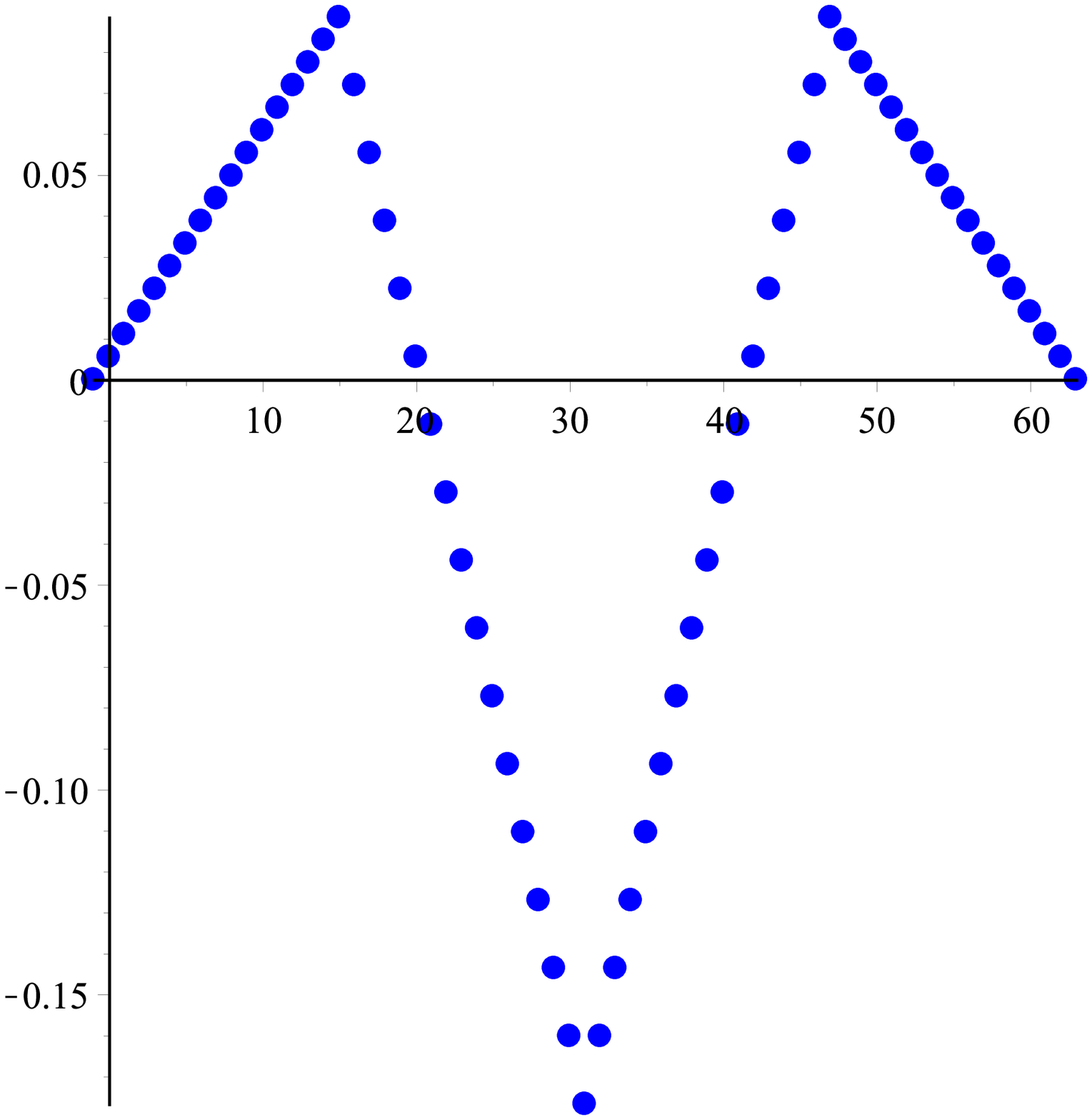}}
\end{tabular}
\captionof{figure}{Plots of $\psi_k^{(1)}$ (left) and $\psi_k^{(2)}$ (right) in
$\ell^2(\mz)$ for $M=10$ and $k=5$ in Example \ref{41228a}.}
\label{fig1}
\end{figure}

As an illustration, consider the functions $\psi_k^{(1)}, \psi_k^{(2)} \in \ell^2(\mz),$ $k =0,\ldots, M-1,$ constructed from the B-splines of
second order $\phi_{k+1,2},$ $k=0,\ldots, M-1,$ as in Proposition \ref{41229c}. Figure \ref{fig1} shows the plots of $\psi_k^{(1)}$ and
$\psi_k^{(2)}$ for $M = 10$ and $k = 5.$
\ep \enx

\bex Let $G=\mr^2,$ with dual group $\hg= \mr^2.$ Consider a $2\times 2$ matrix $A$ with integer entries, eigenvalues outside the unit circle and $| \det A|=2.$
We will show that under a certain technical condition (see \eqref{5101a}) the conditions \eqref{901g} and \eqref{1901n} are satisfied, i.e.,
the B-spline constructions in Proposition \ref{41229c} can be realized. The condition
\eqref{5101a} is satisfied, e.g., for the matrices $A =
\begin{pmatrix} 1 & -1 \\ 1 & 1 \end{pmatrix}$ and $A =
\begin{pmatrix} 0 & 2 \\ -1 & 0 \end{pmatrix}.$
For $k\ge 0,$ consider the lattices
$ \label{4489} \Lambda_k := (A^k)^\sharp \mz^2 = (A^\sharp)^k \mz^2$
in $G= \mr^2,$ where $A^\sharp:=(A^T)^{-1}.$
Observe that as $A$ has integer entries, $A^T \mz^2 \subset \mz^2;$ it follows that
$\mz^2 \subset A^\sharp \mz^2,$ and therefore
$\Lambda_k = (A^\sharp)^k \mz^2 \subset  (A^\sharp)^{k+1} \mz^2= \Lambda_{k+1}, \, \forall k\ge 0,$
i.e., the sequence of lattices $\Lambda_k, \, k\ge 0,$ is nested.
It is easy to check that $| \Lambda_{k+1}/\Lambda_k|=2$ for all $k\ge 0.$
By general theory, the annihilator of the lattice $\lak$ is
$\label{4491} \lak^\bot= \left((A^k)^\sharp\right)^\sharp \mz^2 = A^k \mz^2,$
with associated fundamental domain
$\label{4492} V_k= A^k [0,1)^2.$ For any $k\ge 0,$
$\mu_{\hg}(V_k)= \mu_{{\mr}^2}(A^k [0,1)^2)= | \det A^k| =2^k.$
Let $\{0, \nu\}$ be a full collection of coset representatives of $\mz^2/A\mz^2.$ Then
\bes \mz^2= A \mz^2 \cup \left( \nu + A \mz^2\right), \hspace{1cm} \,
 A \mz^2 \cap \left( \nu + A \mz^2\right)= \emptyset,\ens which implies that
\bes A^k\mz^2= A^{k+1} \mz^2 \cup \left( A^k\nu + A^{k+1} \mz^2\right), \hspace{1cm} \,
 A^{k+1} \mz^2 \cap \left( A^k\nu + A^{k+1} \mz^2\right)= \emptyset.\ens
Thus,
\bes
\lak^\bot= \lako^\bot \cup \left( \nu_k + \lako^\bot\right),
\hspace{1cm} \, \lako^\bot \cap \left( \nu_k + \lako^\bot\right)= \emptyset,\ens
where $\nu_k= A^k \nu \in \lak^\bot \setminus \lako^\bot.$

We will now examine the existence of fundamental domains $Q_k$ for the lattices
$\Lambda_k$ such that \eqref{1901n} is satisfied.
By Theorem 2 in the paper \cite{GM} by Gr\"ochenig and Madych, there exists a
relatively compact set
$Q \subset \mr^2$ and $\eta\in \mr^2$ such that
\bee  \label{4495} Q \cup (\eta + Q) & = &  A^TQ,   \\ \label{44101}
\bigcup_{n\in \mz^2} (n+Q) & = & \mr^2, \\
\label{44104} Q \cap (\eta +Q) & = & \emptyset,  \ene
up to a set of measure zero.
The set $Q$ is generated iteratively by the algorithm
\bes Q^{(r+1)} = A^\sharp Q^{(r)}  \cup A^\sharp (\eta + Q^{(r)}), \, r \geq 0.\ens Also, the  condition  \eqref{4495} is equivalent to
\bee \label{44103} Q= A^\sharp Q\cup A^\sharp(\eta + Q),\ene i.e., $Q$ is self-similar in the sense of the affine transformation $A^\sharp.$

We will now assume that
\begin{equation}
\label{5101a}\mu_{{\mr}^2}(Q)=1;
\end{equation}
Theorem 3 in \cite{GM} provides various equivalent conditions for this condition
to hold. Then Lemma 1 or Theorem 3 in \cite{GM} implies that
$Q \cap (n + Q) = \emptyset, \, \forall n\in \mz^2 \setminus \{0\},$
or, equivalently
\bee \label{44102} (n_1+Q) \cap (n_2+Q)= \emptyset, \, \forall n_1\neq n_2,\, n_1,n_2\in \mz^2,\ene
again up to a set of measure zero.

Applying $(A^\sharp)^k$ to \eqref{44101}, \eqref{44102}, \eqref{44103},
and $(A^\sharp)^{k+1}$  to \eqref{44104} immediately yields that
with $Q_k:= (A^\sharp)^k Q$ and $\eta_k:= (A^\sharp)^{k+1}\eta,$
\bee \label{44105}  \bigcup_{\omega \in \Lambda_k} (\omega + Q_k) & = & \mr^2,
\hspace{.4cm}   (\omega_1 + Q_k)  \cap (\omega_2 + Q_k)  =  \emptyset \, \, \mbox{for} \, \,
\omega_1, \omega_2\in \lak, \omega_1\neq \omega_2, \\
\label{44107}  Q_k & = & Q_{k+1} \cup (\eta_k + Q_{k+1}), \hspace{.4cm}  Q_{k+1} \cap (\eta_k + Q_{k+1})  =  \emptyset.\ene
By \eqref{44105}, the set $Q_k$ is a fundamental domain associated with the lattice $\Lambda_k.$ Also, \eqref{44107} shows that \eqref{1901n} holds.

We now verify the conditions \eqref{2001b} and \eqref{901g}. Let $S\subset \mr^2$ be compact and consider $0< \delta <1.$ Let $\ga \in S,$ and $x\in Q_k;$ we can write $x= (A^\sharp)^kq$ for some $q\in Q.$ Then
\bee \label{41231a} | (-x, \ga)-1| = | e^{-2\pi i \ga \cdot (A^\sharp)^k q}-1|.\ene
Now, using that $S$ and $Q$ are (relatively) compact, and thus bounded, there is a
constant $C>0$ such that, regardless of the choice of $\ga \in S$ and $x\in Q_k,$
\bee \label{41231b} | -2\pi i \ga \cdot (A^\sharp)^k q| \le 2\pi \|\ga\|_2 \, \|(A^\sharp)^k q\|_2
\le 2\pi \|\ga\|_2 \, \|(A^\sharp)^k\|_2 \, \| q\|_2
\le C \, \|(A^\sharp)^k\|_2.\ene
Since $A^T$ has eigenvalues outside the unit circle, there exists a matrix norm
$\| \cdot\|$ such that $\|(A^T)^{-1}\| <1;$ since all norms on a finite-dimensional space are equivalent, this implies that there is a constant $D>0$ such that
\bes  \|(A^\sharp)^k\|_2= \|(A^T)^{-k}\|_2 \le D\, \|(A^T)^{-k}\| \le D\, \|(A^T)^{-1}\|^k;\ens thus, via \eqref{41231b},
\bes  | -2\pi i \ga \cdot (A^\sharp)^k q| \le C D\, \|(A^T)^{-1}\|^k,\, \forall
\ga \in S, x\in Q_k.\ens
Since $\|(A^T)^{-1}\|^k\to 0$ as $k\to \infty,$ it now follows from \eqref{41231a}
that for any given $\delta\in (0,1),$ there exists $k\in I$ such that the inequality \eqref{2001b} holds for all $\ga \in S$ and $x\in Q_k.$

In order to establish \eqref{901g}, for any $k \geq 0, \, n_1, n_2 \in \mz^2,$ consider $\ga \in
 (A^k n_1+S) \cap (A^k n_2 + S).$ Writing
$ \ga = A^k n_1+s_1$ and $\ga = A^k n_2 + s_2,$ where $s_1, s_2 \in S,$ we have  $A^k(n_1-n_2)= s_2-s_1.$ Thus, using that $S$ is
compact, there is a constant $C_1>0$ (independent of the choice of $k$) such that
$ \| A^k(n_1-n_2)\|_2 \le C_1.$
Then
\bes  \| n_1 - n_2 \|_2 = \|(A^k)^{-1} A^k(n_1-n_2)\|_2
& \le & \|(A^k)^{-1}\|_2 \cdot \| A^k(n_1-n_2)\|_2
\le C_1\|(A^k)^{-1}\|_2  \\ & \le &  C_1 \|A^{-1}\|_2^k.\ens
But the eigenvalues of $A$ lie outside the unit circle, so
$\|A^{-1}\|_2^k\to 0$ as $k\to \infty;$ thus for sufficiently large values of $k,$ we have  $n_1=n_2$, which finally proves \eqref{901g}.
\ep \enx

\section{Frames generated by characteristic functions} \label{471225a}
In the entire section, we will use the setup and notations presented in Section \ref{2101a}. In particular,  $\Lambda_k, k\in I,$ is a nested sequence of lattices in $G.$
We consider a nested sequence of subsets $\Omega_k, \, k\in I,$ of $\hg$ such that for $k \in I,$
\bee \label{4454c} \Omega_k \subseteq V_k \ene
for some fundamental domain $V_k$ associated with $\Lambda_k^\bot.$ We
will further assume that for every compact set $S$ in $\hg,$ there exists $K \in I$ such that $S\subseteq \Omega_K.$
Now, for $k\in I,$ let
\bee \label{4457} \Phi_k(\ga):= \muhg(V_k)^{-1/2} \chi_{\Omega_k}(\ga), \, \ga \in \hg.\ene
We will show that with this setup, we can always satisfy the UEP conditions, regardless
of the underlying LCA group.
As in Section \ref{2101a}, with the set $V_k$ and $\nu_{k, \ell}$ defined as in \eqref{1001d}, the fundamental domain
$V_{k+1}^\prime =  \bigcup_{\ell=1}^{d_k} ( \nu_{k,\ell}+V_k)$ in \eqref{172c} associated with $\Lambda_{k+1}^\bot$ will be useful in our derivations.
Note that in contrast with the B-spline generated case in
Section \ref{1901d}, we assume neither  the condition \eqref{41229a} on the lattices, nor  the splitting \eqref{1901n} of the fundamental domains $Q_k$ associated with
the lattices $\Lambda_k.$

We will first show that  in the above setup, the refinement equation
\eqref{1001a} and the technical conditions \eqref{901g} and \eqref{901a} are always satisfied.
\bl \label{41227a} Let  $\{ \Lambda_k\}_{k\in I}$ in $G$ be a nested sequence of lattices in $G,$ and choose the sets $V_k$ and $\Omega_k$ as above.
Then the functions $\Phi_k,$ $k \in I,$ in \eqref{4457} satisfy the refinement equation
\eqref{1001a} whenever $H_{k+1}\in L^\infty(V_{k+1})$ is defined by \bee \label{4459}\hko(\ga):= \begin{cases} \sqrt{d_k}, \, &\mbox{if} \, \ga \in \ok,\\
0 , \, &\mbox{if} \, \ga \in \vko^\prime \setminus \ok. \end{cases}\ene
Furthermore, the conditions \eqref{901g} and \eqref{901a} are satisfied. \el

\bp  For $k \in I,$ if $\ga \in \Omega_k,$ then $\Phi_k(\ga)= \muhg(V_k)^{-1/2}$
and $\Phi_{k+1}(\ga)= \muhg(V_{k+1})^{-1/2},$ i.e., \eqref{1001a} holds if we
let
\bes H_{k+1}(\ga)= \frac{\Phi_k(\ga)}{\Phi_{k+1}(\ga)}=\frac{ \muhg(V_k)^{-1/2}}{ \muhg(V_{k+1})^{-1/2}}= \left(\frac{s(\Lambda_k^\bot)}{s(\Lambda_{k+1}^\bot)}\right)^{-1/2}= \left(\frac{s(\Lambda_{k+1})}{s(\Lambda_k)}\right)^{-1/2}
=\sqrt{d_k}, \, \ga \in \Omega_k.\ens
If $\ga \in V_{k+1}^\prime \setminus \ok,$ then $\Phi_k(\ga) = 0$ and \eqref{1001a} is satisfied when
we take $H_{k+1}(\ga)=0,$ regardless of the value of $\Phi_{k+1}(\ga).$
We can therefore
extend the function $\hko$ in \eqref{4459} to a periodic function, and
\eqref{1001a} holds.

We will now check the technical conditions \eqref{901g} and \eqref{901a}. First,
given a compact set $S$ in $\hg,$ it follows from the inclusions $S \subseteq \Omega_K \subseteq V_K$ that
for any $\omega, \omega^\prime \in \Lambda_K^\bot,$ $\omega \neq \omega^\prime,$
$(\omega + S) \cap (\omega^\prime + S) \subseteq
(\omega + V_K) \cap (\omega^\prime + V_K).$ Since $\mu_{\hg}\big( (\omega + V_K) \cap (\omega^\prime + V_K) \big) =0,$
this gives \eqref{901g}. Next,
the definition of $\Phi_k$ and
the assumptions of the sets $\Omega_k$ imply that for $k\ge K$ and $\ga \in S,$
$  | \mu_{\hg}(V_k) |\Phi_k(\ga) |^2 -1 |=0,$ i.e., \eqref{901a} is satisfied.
\ep



With the considered setup we can apply the UEP as soon as we have
defined functions $G_{k+1}^{(m)}\in L^\infty(\vko), \,
m=1, \dots, \rho_k,$ in such a way that the matrix $P_k$ in \eqref{1001g}, with $\rho_k$ equals $d_k$ or $d_k -1$, satisfies
\eqref{172e}.  We note that it is enough to define $G_{k+1}^{(m)}(\ga + \nu_{k, \ell})$
for $\ga \in \vk, m=1, \dots, \rho_k, \ell=1, \dots, d_k;$ then, if for
$\omega\in \Lambda_{k+1}^\bot$ we define
\bee \label{4464} G_{k+1}^{(m)}(\ga + \nu_{k, \ell}+\omega):= G_{k+1}^{(m)}(\ga + \nu_{k, \ell}),\ene the entries of $P_k$ are defined as periodic functions in $L^\infty(\vko).$

In order to consider the matrix extension problem \eqref{172e}, we need an explicit expression for the first
row of the matrix $P_k$, which follows immediately from Lemma \ref{41227a}.

\bl \label{41208a} Let $\hko$ be defined by \eqref{4459}. Then for $\ga \in \vk,$
\bes
(\hko(\ga+ \nu_{k,1}), \dots, \hko(\ga+ \nu_{k,d_k}))=(\hko(\ga),0, \dots,0).\ens  \el

We will now construct explicitly given tight frames based on the
functions $\Phi_k$ in \eqref{4457}. These tight frames are ``bandlimited'' in the sense that they are compactly supported on
the dual group $\hg.$ We first consider the case where $\Omega_k$ is
chosen as a proper subset of $V_k.$

\bpr \label{41208g} Let  $\{ \Lambda_k\}_{k\in I}$ be a nested sequence of lattices in $G,$ and choose the sets $V_k$ and $\Omega_k$ as above; assume further
that
$\Omega_k$ is a proper subset of $V_k$ for all $k\in I.$
For $k \in I,$ consider the function $\Phi_k$ in \eqref{4457} with the associated filter
$H_{k+1}\in L^\infty(V_{k+1})$ in \eqref{4459}. Define the functions
$G_{k+1}^{(m)}, m=1, \dots, d_k-1,$ by
\bes G_{k+1}^{(m)}(\ga + \nu_{k, \ell})  :=  \begin{cases} \sqrt{d_k} \delta_{m+1, \ell}, \, & \mbox{if} \,\ga \in \Omega_k, \\
\sqrt{d_k} \delta_{m, \ell}, \,   & \mbox{if} \,\ga \in \vk \setminus \ok,
\end{cases} \ens and the function $G_{k+1}^{(d_k)}$ by
\bes G_{k+1}^{(d_k)}(\ga + \nu_{k, \ell}):= \begin{cases} 0, \, & \mbox{if} \,\ga \in \Omega_k, \\  \sqrt{d_k} \delta_{d_k, \ell}, & \mbox{if} \,\ga \in \vk \setminus \ok,
   \end{cases} \ens
where $\ell=1, \dots, d_k.$ Then the collection of functions
\bes  \{\ml \Phi_{k_0}\}_{\lambda\in \Lambda_{k_0}} \bigcup \,
\{\ml \Psi_k^{(m)}\}_{k\ge k_0, \lambda\in \Lambda_{k}, m=1, \dots, d_k}\ens forms a tight frame for $\lthg$ with frame bound $1,$ and the functions  $\Psi_k^{(m)}$ are compactly
supported. \epr

\bp Recall that
with $\hko$  defined by \eqref{4459}, Lemma \ref{41208a} determines the first row of
the matrix $P_k(\ga)$ for $\ga \in V_k.$
Then by considering $\ga \in \Omega_k,$ we obtain
\bes P_k(\ga)= \begin{pmatrix}  \sqrt{d_k} I_{d_k} \\ {\bf 0} \end{pmatrix}, \, \ga \in \ok,\ens where ${\bf 0}$ denotes the zero vector in $\mc^{d_k}.$ Thus,
$P_k(\ga)^* P_k(\ga)=d_k I_{d_k}, \, \ga \in \ok.$ Next, by considering
$\ga \in \vk \setminus \ok,$ we have
 \bes P_k(\ga)= \begin{pmatrix}  {\bf 0} \\ \sqrt{d_k} I_{d_k} \end{pmatrix}, \, \ga \in \vk \setminus \ok,\ens which
implies that $P_k(\ga)^* P_k(\ga)=d_k I_{d_k}, \, \ga \in \vk \setminus \ok.$ Altogether,
this verifies \eqref{172e}. Thus, by Lemma \ref{41227a}  and Theorem \ref{1001h}, the collection of functions
\bes  \{\ml \Phi_{k_0}\}_{\lambda\in \Lambda_{k_0}} \bigcup \,
\{\ml \Psi_k^{(m)}\}_{k\ge k_0, \lambda\in \Lambda_{k}, m=1, \dots, d_k}\ens forms a tight frame for $\lthg$ with frame bound 1.
Since
\bes \Psi_k^{(m)} (\ga)= G_{k+1}^{(m)}(\ga) \Phi_{k+1}(\ga)= \muhg(\vko)^{-1/2}
G_{k+1}^{(m)}(\ga) \chi_{\oko}(\ga),\ens
we see that $\supp\, \Psi_k^{(m)}\subseteq \oko \subset \vko,$ i.e., $\Psi_k^{(m)}$
is compactly supported.
\ep

Note that the construction in Proposition \ref{41208g} forms only one way of extending the first row of the matrix $P_k(\ga), \, \ga \in \vk$
in order to obtain  \eqref{172e}.  Indeed, observe that
\bes \sum_{\ell=1}^{d_k} | \hko (\ga + \nu_{k, \ell})|^2 = \begin{cases}
d_k, \ &\mbox{if} \, \ga \in \ok, \\
0, \ &\mbox{if} \, \ga \in \vk \setminus \ok,\end{cases} \ens
i.e.,
$\sum_{\ell=1}^{d_k} | \hko (\ga + \nu_{k, \ell})|^2 \le d_k, \, \forall \ga \in V_k.$ Thus, for any $\rho_k \ge d_k,$ by using the theory of
Householder matrices one can construct explicitly $(\rho_k+1) \times d_k$ matrices $P_k(\ga)$ that satisfy
\eqref{172e}. Details can be found in Proposition 4.1 of \cite{GT2} which deals with the group $\mt^s$, but its ideas can be readily adapted to the current setting of general LCA groups.

We will now consider the choice $\Omega_k:=V_k,$ which leads to an analogue on LCA groups of the Shannon orthonormal basis for $L^2(\mr)$.

\bpr \label{41208c} Given a nested sequence of lattices $\{ \Lambda_k\}_{k\in I}$ in $G,$ let $V_k$ be a fundamental domain associated  with the lattice
$\Lambda_k^\bot$ in $\hg,$ and suppose that the sets $V_k,$ $k \in I,$ are nested.
Assume that
for every compact set $S$ in $\hg,$ there exists $K \in I$ such that $S\subseteq V_K.$
For $k\in I,$ let
\bee \label{4457a} \Phi_k(\ga):= \muhg(V_k)^{-1/2} \chi_{V_k}(\ga), \, \ga \in \hg,\ene
with the associated filter $H_{k+1}$ as in \eqref{4459} $($with $\Omega_k=V_k).$
For $m=1, \dots, d_k-1,$ define the function $G_{k+1}^{(m)}$ by
\bee \label{4475} G_{k+1}^{(m)}(\ga + \nu_{k, \ell}):= \sqrt{d_k}\delta_{m+1, \ell}, \, \ga \in V_k,\ene
for $\ell=1, \dots, d_k.$ Then the collection of functions
\bee \label{41230ff}  \{\ml \Phi_{k_0}\}_{\lambda\in \Lambda_{k_0}} \bigcup \,
\{\ml \Psi_k^{(m)}\}_{k\ge k_0, \lambda\in \Lambda_{k}, m=1, \dots, d_k-1}\ene forms an orthonormal basis for $\lthg.$
\epr

\bp Using Lemma \ref{41208a}, we may simply require that
$\label{4473} P_k(\ga):= \sqrt{d_k}I_{d_k}, \, \ga \in \vk,$
which leads to the choice of $G_{k+1}^{(m)}$ in \eqref{4475}.
 By periodicity as in \eqref{4464},
the matrix $P_k$ in \eqref{1001g} satisfies
\eqref{172e}. Thus, by Lemma \ref{41227a} and Theorem \ref{1001h},  the functions
 in \eqref{41230ff} form a tight frame for $\lthg$ with frame bound 1. Furthermore, by \eqref{4457a},
$\| \Phi_{k_0}\|=1.$ Also, for $k\in I, m=1, \dots, d_k-1,$
$\Psi_k^{(m)}(\ga)= G_{k+1}^{(m)}(\ga) \Phi_{k+1}(\ga)=\muhg(V_{k+1})^{-1/2} G_{k+1}^{(m)}(\ga) \chi_{V_{k+1}}(\ga);$ thus
\bes \|\Psi_k^{(m)}\|^2 = \muhg(V_{k+1})^{-1}\int_{\vko} | G_{k+1}^{(m)}(\ga)|^2\, d\ga.\ens
As $V_k$ and $V_{k+1}^\prime$ are fundamental domains associated with $\Lambda_k^\bot$ and $\Lambda_{k+1}^\bot$ respectively, we
may express $V_{k+1}$ as the disjoint union $V_{k+1} = \bigcup_{\omega \in \Lambda_{k+1}^\bot}
\bigcup_{\ell = 1}^{d_k} ( \omega + \nu_{k, \ell} + V_k) \cap V_{k+1}.$ Using this decomposition of $V_{k+1},$ it follows that
\bes \|\Psi_k^{(m)}\|^2   =  \muhg(V_{k+1})^{-1} \sum_{\omega \in \Lambda_{k+1}^\bot}
\sum_{\ell=1}^{d_k}\int_{(\omega + \vk) \cap (-\nu_{k, \ell} + V_{k+1})} \hspace*{-0.8ex}
| G_{k+1}^{(m)}(\ga+\nu_{k, \ell})|^2\, d\ga,\ens
and therefore by \eqref{4475},
\bes \|\Psi_k^{(m)}\|^2 & = & \muhg(V_{k+1})^{-1} d_k \sum_{\omega \in \Lambda_{k+1}^\bot} \mu_{\widehat G}
( (\omega + \vk) \cap (-\nu_{k, m+1} + V_{k+1}) ) \\
& = & \muhg(V_{k+1})^{-1} d_k \sum_{\omega \in \Lambda_{k+1}^\bot} \mu_{\widehat G}
( \vk \cap (-\omega -\nu_{k, m+1} + V_{k+1}) ) \\ & = &
\muhg(V_{k+1})^{-1} d_k \, \muhg(\vk)=1.\ens Hence, the collection
actually forms an orthonormal basis for $\lthg,$ as claimed. \ep

We will now apply the results in this section to obtain explicit constructions of frames for several concrete LCA groups.
In particular we will be able to verify all the technical assumptions made in the section.

\bex \label{50104f} Let $M\in \mn,$ and consider $G= \mz_{2^M}=\{0, \dots, 2^M-1\}.$ Then
$\hg= G=\mz_{2^M},$ and $L^2(G)= L^2(\hg)= {\cal S}(2^M),$ the space of complex sequences indexed by $\mz$ and with period $2^M.$
For $k=0, \dots, M,$ consider the lattice
$\label{44116} \lak:= 2^{M-k}\mz_{2^k}$ in $G;$ we note that
$\Lambda_0= \{0\}$ and $\Lambda_M= G.$ Using that
$ 2\mz_{2^k} =\{0, 2, \dots, 2^{k+1}-2 \}\subset \{0, 1, \dots, 2^{k+1}-1\}=\mz_{2^{k+1}},$
$\lak =  2^{M-k-1} 2\mz_{2^k} \subset 2^{M-k-1}\mz_{2^{k+1}}= \lako,$ i.e., \eqref{901k} holds. It is well known that
$\label{44117} \lak^\bot= 2^k \mz_{2^{M-k}}, \, k=0, \dots, M.$ Direct calculation verifies that for $k=0, \dots, M-1,$
$  \lak^\bot = \lako^\bot \cup (2^k + \lako^\bot)$ and $\lako^\bot \cap (2^k + \lako^\bot)= \emptyset,$ i.e.,
\eqref{1001d} holds with $\nu_{k,0}=0, \nu_{k,1}=2^k.$ Also, $V_0=\{0\}$ is a fundamental domain associated with the lattice $\Lambda_0^\bot,$ and for $k=0, \dots, M-1,$ the set
$ \label{44119} V_{k+1}:= \mz_{2^{k+1}}$ is a fundamental domain for $\lako^\bot.$

In order to define appropriate sets $\Omega_k, k=0,\dots,M,$ let
$\{L_k\}_{k=0, \dots, M}$ be an increasing sequence of nonnegative integers satisfying that
$$ L_M=2^M-1, \hspace{.4cm} L_k \le 2^k-1, \, k=0,\dots, M-1.$$
Now, for $k=0, \dots, M,$ let
$$ \Omega_k:= \{0, \dots, L_k\}.$$
Since $L_k\le 2^k-1$ and $V_k= \mz_{2^k}$, it follows that
$\Omega_k \subseteq V_k,$ i.e., \eqref{4454c} holds. Also, since $\{L_k\}_{k=0,\dots,M}$ is increasing, we clearly have $\Omega_k \subseteq \Omega_{k+1}.$ By the choice of $L_M=2^M-1,$
for any compact set $S$ in $\hg$ we have $S \subseteq \Omega_M= \hg.$ Thus, we have verified all the assumptions for application of the UEP. \ep \enx

\bex \label{50106a} Let $G= \mt;$ then $\hg= \mz.$ Let $\{M_k\}_{k\ge 0}$ denote a sequence of integers $M_k\ge 2,$ and assume for convenience that $M_0$ is an even number. Define
$ \label{44123} N_k:= \prod_{\ell=0}^k M_\ell, \, k\ge 0,$ which are even integers
greater than or equal to $2.$ For $k\ge 0,$ consider the lattice
$\label{44124} \lak:= \frac1{N_k} \mz_{N_k}$ in $G;$ then
$\lak = \frac1{\nk \mko}\, \mko \mz_{\nk} \subset   \frac1{\nk \mko}\,  \mz_{\nk \mko}
= \lako,$
i.e., the lattices are nested. Also, for $k\ge 0,$
$ \label{44125} \lak^\bot = \nk \mz,$ and direct verification yields that
$\lak^\bot  =  \bigcup_{\ell=1}^{\mko} \left((\ell -1)\nk + \lako^\bot\right)$ and \bes \left((\ell -1)\nk + \lako^\bot \right) \cap
\left((\ell^\prime -1)\nk + \lako^\bot \right) = \emptyset, \,  \ell\neq \ell^\prime, \, \ell, \ell^\prime=1, \dots, \mko,\ens
where $(\ell -1)\nk \in \lak^\bot \setminus \lako^\bot, \, \ell=2, \dots,
\mko.$

For $k\ge 0,$ the set
$ \label{44126} V_k:= \{ -\frac{\nk}{2} ,  - \frac{\nk}{2}  +1, \dots,    \frac{\nk}{2}  -1   \}$
is a fundamental domain associated with $\lak^\bot.$ Now, let $\{L_k\}_{k\ge 0}$ be a
(strictly) increasing sequence of nonnegative integers such that
$L_k \le \frac{\nk}{2}-1, \, k\ge 0,$ and define
$$ \Omega_k:= \{-L_k, \dots, L_k\}, \, k\ge 0.$$
Then $\Omega_k \subset V_k,$ and also $\Omega_k \subset \Omega_{k+1}.$ Finally,
if $S\subset \hg= \mz$ is compact, i.e., $S\subset \{-R, \dots, R\}$ for some $R>0,$
taking $K \geq 0$ such that $L_K>R$ implies that $S\subset \Omega_K;$ thus, all the assumptions for application of the UEP are verified. \ep \enx

In our final example, we consider $s\times s$ diagonal scaling matrices and construct nonstationary wavelet frames for $\ltrs$ (separable as well as nonseparable).

\bex \label{50106e} Let $G= \mrs;$ then $\hg= \mrs.$  For $r=1, \dots, s,$ let
$\{M_{k,r}\}_{k\ge 0}$ be a sequence of integers $M_{k,r}\ge 2.$ For $k\ge 0,$ let
$ \label{44130} N_{k,r}:= \prod_{\ell=0}^k M_{\ell,r},$ and consider the
$s\times s$ matrix
$ \label{44131} A_k:= \mbox{diag} \, (N_{k,1}, \dots, N_{k,s}).$
We will now consider the lattice
$ \label{44132} \lak:= A_k^\sharp \mzs= \frac1{N_{k,1}}\, \mz \times  \cdots
\times \frac1{N_{k,s}}\, \mz;$ introducing the short notation
$ \prod_{r=1}^s B_r:= B_1 \times \cdots \times B_s$ for the cartesian product of $s$
sets,
the lattice can be written as
$$ \label{44132a} \lak= \prod_{r=1}^s \left(\frac1{N_{k,r}}\, \mz\right).$$
Now,
\bes \lak =  \prod_{r=1}^s \left(\frac1{N_{k,r} M_{k+1,r}}\, M_{k+1,r} \mz\right)
=  \prod_{r=1}^s \left(\frac1{N_{k+1,r} }\, M_{k+1,r} \mz\right)
\subset \prod_{r=1}^s \left(\frac1{N_{k+1,r} }\,  \mz\right)= \lako,\ens
i.e., the sets $\lak,$ $k \geq 0,$ are nested. Also, for $k\ge 0,$
$ \label{44133} \lak^\bot = \prod_{r=1}^s (N_{k,r}\, \mz).$
Direct verification yields that
$ \lak^\bot= \bigcup_{\ell=1}^{d_k} \left( \nu_{k,\ell} + \lako^\bot\right),$
where $d_k= \prod_{r=1}^s M_{k+1,r},$ and
$$ \label{44134} \{\nu_{k,\ell}\}_{\ell=1}^{d_k}= \prod_{r=1}^s \left(N_{k,r}
\{0, 1, \dots, M_{k+1,r}-1\}\right)$$
and $\nu_{k, \ell} \in \lak^\bot \setminus \lako^\bot, \, \ell=2, \dots, d_k;$
finally,
\bes \left( \nu_{k, \ell}+ \lako^\bot \right) \cap  \left( \nu_{k, \ell^\prime}+ \lako^\bot \right)= \emptyset, \,
 \ell \neq \ell^\prime, \, \ell, \ell^\prime=1, \dots, d_k.\ens For $k\ge 0,$ the set
$\label{44135} V_k:= \prod_{r=1}^s \big[ - \frac{N_{k,r}}{2}, \frac{N_{k,r}}{2}\big)$ is a fundamental domain associated with $\lak^\bot.$

We shall now identify sets $\Omega_k, \, k\ge 0,$ satisfying our technical conditions.
For $r=1, \dots, s,$ let $\{L_{k,r}\}_{k\ge 0}$ be an increasing sequence of nonnegative numbers (not necessarily integers) satisfying that
$ L_{k,r}\le \frac{N_{k,r}}{2}, \, \forall k\ge 0,$
and
$ \lim_{k\to \infty} L_{k,r}= \infty.$
Then, for $k\ge 0,$ we can take
$ \label{44138} \Omega_k:= \prod_{r=1}^s [ - L_{k,r}, L_{k,r});$  we leave the easy proof
to the reader. This yields separable wavelet frames.

For the construction of nonseparable wavelet frames, let $\{L_k\}_{k\ge 0}$ be an increasing sequence of nonnegative numbers
(not necessarily integers) satisfying that
$$ \label{44139} L_k < \frac12 \min(N_{k,1}, \dots N_{k,s}), \, \forall k\ge 0,
\hspace{.2cm} \mbox{and} \hspace{.3cm} \lim_{k\to \infty} L_k= \infty.$$
For $k\ge 0,$ let
$$ \label{44141} \Omega_k:= \{ \ga \in \mrs \, \big| \, \|\ga\|_2 \le L_k\}.$$
For $\ga=(\ga_1, \dots, \ga_s)\in \Omega_k,$ each coordinate $\ga_r$ satisfies that
$ | \ga_r| \le L_k < \frac{N_{k,r}}{2},$
i.e., $\ga \in V_k;$ thus $\Omega_k \subset V_k,$ i.e., \eqref{4454c} holds.
It is also clear that $\Omega_k \subseteq \Omega_{k+1}$ for all $k\ge 0.$ Finally, let
$S\subset \hg = \mrs$ be a compact set. Then, choosing $R>0$ such that $S$ is contained in the ball around 0 with radius $R$ and choosing $K\ge 0$ such that $L_K >R,$
we clearly have that $S\subset \Omega_K;$ thus, again, all the technical assumptions
for application of the UEP are satisfied.
\ep \enx

\vspace{.1in}\noindent{\bf Acknowledgment:} Ole Christensen would like to thank the National University of Singapore for its warm hospitality during one month stays in December 2013 -- January 2014, December 2014 -- January 2015, and December 2015 -- January 2016.

\begin{tabbing}
text-text-text-text-text-text-text-text-text-text \= text \kill \\
Ole Christensen \> Say Song Goh \\
Technical University of Denmark \> Department of Mathematics \\
DTU Compute \> National University of Singapore \\
Building 303, 2800 Lyngby \> 10 Kent Ridge Crescent \\
Denmark \> Singapore 119260, Republic of Singapore \\
Email: ochr@dtu.dk \> Email: matgohss@nus.edu.sg
\end{tabbing}


\begin{thebibliography}{10}

\frenchspacing

\bibitem{Skopina} S.\ Albeveiro, S.\ Evdokimov, M.\ Skopina,
$p$-adic multiresolution analysis and wavelet frames, J.\ Fourier Anal.\
Appl.\ {16} (2010) 693--714.

\bibitem{benedetto} J.\ Benedetto, R.L.\ Benedetto,  A wavelet theory for
local fields and related groups, J.\ Geom.\ Anal.\ {14} (2004)
423--456.

\bibitem{BowRos2014}
M.\ Bownik, K.A.\ Ross,
The structure of translation-invariant spaces on locally compact
  abelian groups, J.\ Fourier Anal.\ Appl.\ 21 (2015) 849--884.


\bibitem{CP} C.\ Cabrelli, V.\ Paternostro, {Shift-invariant spaces on LCA groups,} J.\ Funct.\ Anal.\ {258} (2010) 2034--2059.

\bibitem{CBN} O.\ Christensen, An Introduction to Frames and Riesz Bases, 2nd
expanded ed.,
Birkh\"auser, 2016.

\bibitem{CG} O.\ Christensen, S.S.\ Goh,
{Fourier-like frames on locally compact abelian groups,}  J.\ Approx.\ Theory 192 (2015) 82--101.


\bibitem{CHS}
C.K.\ Chui, W.\ He, J.\ St\"{o}ckler, Compactly supported tight
and sibling frames with maximum vanishing moments,
Appl.\ Comput.\ Harmon.\ Anal.\
13 (2002) 224--262.

\bibitem{Dah} S.\ Dahlke,
Multiresolution analysis and wavelets on locally compact abelian groups,
in:
P.-J.\ Laurent, A.\ Le M\'{e}haut\'{e}, L.\ Schumaker (Eds.),
Wavelets, Images, and Surface Fittings, AK Peters, 1994, pp.\ 141--156.

\bibitem{DRoSh3}
I.\ Daubechies, B.\ Han, A.\ Ron, Z.\ Shen,
Framelets: MRA-based constructions of wavelet frames,
Appl.\ Comput.\ Harmon.\ Anal.\
14 (2003) 1--46.


\bibitem{FeiKo} H.G.\ Feichtinger, W.\ Kozek, {Quantization of TF lattice-invariant operators on elementary LCA groups,}
    in: H.G.\ Feichtinger, T.\ Strohmer (Eds.),
Gabor Analysis and Algorithms: Theory and Applications, Birkh\"auser,
1998, pp.\ 233--266.



\bibitem{GT2} S.S.\ Goh, K.M.\ Teo,
Extension principles for tight wavelet frames of periodic functions,
Appl.\ Comput.\ Harmon.\ Anal.\ 25 (2008) 168--186.

\bibitem{G6} K.\ Gr\"ochenig,
{Aspects of Gabor analysis on locally compact abelian groups,} in: H.G.\ Feichtinger, T.\ Strohmer (Eds.),
Gabor Analysis and Algorithms: Theory and Applications, Birkh\"auser,
1998, pp.\ 211--231.

\bibitem{GM} K.\ Gr\"ochenig, W.R.\ Madych,
{Multiresolution analysis, Haar bases, and self-similar tilings of $\mr^n,$}
IEEE Trans.\ Inform.\ Theory {38} (1992) 556--568.


\bibitem{HeRo} E.\ Hewitt, K.A.\ Ross, {Abstract Harmonic Analysis,
vol.\ 1 and 2,} Springer, 1963.



\bibitem{LeSi1} M.S.\ Jakobsen, J.\ Lemvig,
{Reproducing formulas for generalized
translation invariant systems on locally compact groups,} Trans.\ Amer.\ Math.\ Soc.\ {368} (2016) 8447--8480.



\bibitem{KK} E.\ Kaniuth, G.\ Kutyniok, {Zeroes of the Zak transform on locally compact abelian groups,} Proc.\ Amer.\ Math.\ Soc.\
{126} (1998) 3561--3569.

\bibitem{KL} G.\ Kutyniok, D.\ Labate, {Theory of reproducing systems on locally compact abelian group,}  Colloq.\ Math.\ {106} (2006) 197--220.


\bibitem{ReSt} H.\ Reiter, J.D.\ Stegeman,
{Classical Harmonic Analysis and Locally Compact Groups,} 2nd ed., Oxford University Press, 2000.



\bibitem{RoSh2} A.\ Ron, Z.\ Shen, Affine systems in
$L_2(\mr^d)$: the analysis of the analysis operator,
J.\ Funct.\ Anal.\ 148 (1997) 408--447.





\bibitem{RoSh7} A.\ Ron, Z.\ Shen, {Generalized shift-invariant systems,} Constr.\ Approx.\ {22} (2005) 1--45.



\bibitem{Ru3} W.\ Rudin, {Fourier Analysis on Groups.} Interscience
Publishers, 1962.


\bibitem{S94} G.\ Steidl, {Wavelets over $\mr, \mz, \mr/N\mz$ and $\mz/N\mz$,}
in: C.K.\ Chui, L.\ Montefusco, L.\ Puccio (Eds.), Wavelets: Theory,
Algorithms and Applications, Academic Press, 1994, pp.\ 155--179.

\bibitem{T} V.M.\ Tikhomirov, {Harmonic tools for approximation and splines on locally compact abelian groups,}
Uspekhi Mat.\ Nauk {49} (1994), 193--194; translated in Russian Math.\ Surveys {49} (1994), 200--201.

\bibitem{WZ} Y.G.\ Wang, X.\ Zhuang, {Tight framelets and fast framelet
transforms on manifolds,} preprint, 2016.

\end{thebibliography}
\end{document}